\documentclass{amsart}

\usepackage{amsmath}
\usepackage{amsfonts}
\usepackage{epsfig}
\usepackage{graphics} 

\usepackage{pgfplots,pgfplotstable}
\usepgfplotslibrary{groupplots}
\usetikzlibrary{shapes,arrows,snakes,calendar,matrix,backgrounds,folding,calc,positioning,spy}
\usepackage{tikz}
\usepackage{booktabs}

\usepackage{siunitx}
\usepackage[utf8]{inputenc}
\usepackage[T1]{fontenc}
\usepackage{graphicx}
\usepackage{amssymb,amsmath}
\usepackage{url}
\usepackage{lastpage}

\usepackage{hyperref}
\definecolor{darkblue}{rgb}{0,0,.7}
\hypersetup{colorlinks=true,linkcolor=darkblue,citecolor=darkblue}

\usepackage[scaled=0.9]{helvet}

\usepackage{eqnarray}

\def\norm#1{\|#1\|}
\def\Norm#1{\vert \!\vert \!\vert#1\vert \!\vert \!\vert}
 \def\Normh#1{\vert \!\vert \!\vert#1\vert \!\vert \!\vert_h}

\newcommand{\vertiii}[1]{{\left\vert\kern-0.25ex\left\vert\kern-0.25ex\left\vert #1 
    \right\vert\kern-0.25ex\right\vert\kern-0.25ex\right\vert}}

\def\enorm#1{\|#1\|_{\mathcal C}} 
\def\hnorm#1{\|#1\|_h} 
\def\anorm#1{\|#1\|_{\mathcal A}} 
\def\snorm#1{\mu^{-1/2}\|#1\|_{0}}

  \newcommand\NEW{\textrm{RAFW}}
  \newcommand\Nel{N}
  \def\korn{C_\Omega} 

   \newcommand\Hdnull{ H_0(\div\!\!:\!\Omega)} 
      \newcommand\Hdg{ H_g(\div\!\!: \!\Omega) } 
      
       \newcommand\Vtwo{[L^2(\Omega)]^\dim } 
  
  \newcommand\Hdiv{ H(\mathrm{div}\!: \!\Omega)}
  
      \newcommand\RTk{RT_{k-1}(K)}
          \newcommand\RTkH{\widetilde {RT}_{k-1}(K)}
             \newcommand\Rl{[P_l(K)]^\dd_{{\rm skw}} }

      \newcommand\Mk{ M_{k+1}(K) }
    \newcommand\Ml{ M_{l+1}(K) }   
      \newcommand\Mtwo{ M_2(K) }

        \newcommand\dd{{d\times d} }
  
  \newcommand\E{{E(\Omega)}}

      \newcommand\Sh{S_h}
   \newcommand\Vh{V_h}

    \newcommand\HH{ [L^2(\Omega)]^{\dim \times \dim}   } 
       \newcommand\SK{ [L^2(\Omega)]^{\dim \times \dim} _{{\rm skw}}  } 
             \newcommand\SY{ [L^2(\Omega)]^{\dim \times \dim} _{{\rm sym}}  } 
   \newcommand\VV{[H_D^1(\Omega)]^\dim}

\newcommand\etth{{h}}
   
   \newcommand\bfs{\sigma}
  \newcommand\bft{\tau}
  
  \newcommand\st{\tilde\sigma}
  \newcommand\zt{z}
  \newcommand\rt{\xi}

 \newcommand\ph{u_h}
  \newcommand\uoh{ u_h^{a}}
  \newcommand\Ioh{I_h^a}
  
\newcommand\Ch{{\mathcal C}_h}

\newcommand\B{{\mathcal B}}
\newcommand\M{{\mathcal M}}

\renewcommand\H{H({\rm div}\!:\!\Omega)}

\newcommand\comp {\mathcal {C}} 
\newcommand\elas {{\mathcal A}} 
\newcommand\tr{{\rm tr}}

\def\BF0{0}

\newcommand{\eps}{\varepsilon}

\renewcommand{\div}{\mathrm{div}\,}
\newcommand{\jump}[1]{ [\! [#1 ]\!]}
\def\real{\mathbb{R}}

\newcommand\Nh{N_{h}}

 \newcommand\curl{{\rm curl}\,}

  \def\dim{d} 
  \def\D{D} 
  
  \newtheorem{theorem}{Theorem}
\newtheorem{lemma}{Lemma}

\newtheorem{remark}{Remark}

\numberwithin{equation}{section}

\usepackage{stmaryrd}

\usepackage{todonotes}

\title[
Weakly Symmetric Mixed  Elements for Elasticity]{Analysis of 
Weakly Symmetric Mixed Finite Elements for Elasticity 
 }

\author[P.~L.~Lederer]{Philip L. Lederer}
\address{Department of Mathematics and Systems Analysis, Aalto University,
Otakaari 1, Espoo, Finland}
\email{philip.lederer@aalto.fi}

\author[R.~Stenberg]{Rolf  Stenberg}
\address{Department of Mathematics and Systems Analysis, Aalto
University, Otakaari 1, Espoo, Finland}
\email{rolf.stenberg@aalto.fi}

\thanks{This work was supported by the Academy of Finland (Decision 324611).}

\begin{document}
\begin{abstract}
We consider mixed finite element methods for linear elasticity where
the symmetry of the stress tensor is weakly enforced. Both an a priori
and a posteriori error analysis are given for several known families
of methods that are uniformly valid in the incompressible limit.   A posteriori estimates
are derived for both the compressible and incompressible cases.  The results are verified
by numerical examples.
\end{abstract}
 

\maketitle
 
\section{Introduction} 

The purpose of this paper is to give a complete error analysis of mixed finite element methods for linear elasticity where the symmetry of the stress tensor is weakly enforced. In an earlier paper \cite{LS1} we have analyzed methods with an exact enforcement of symmetry.

The study of mixed methods with the stress tensor and displacement as unknowns goes back to the pioneering work of B.
Fraeijs de Veubeke \cite{FdV}. When the stress tensor is assumed to be exactly symmetric, this leads to methods that are quite complicated and cumbersome to implement. This was
also known to Fraeijs de Veubeke, and it lead him to introduce methods
for which the symmetry of the stress tensor is enforced in a weak form
\cite{FdV-W,FdV2,FdV3}. The Lagrange multiplier used has the physical
meaning of the infinitesimal rotation and the corresponding Euler
equation is that of rotational equilibrium.

The mathematical analysis of this class of methods was initiated by
Amara and Thomas \cite{AT} who considered  hybrid methods, and Arnold,
Brezzi and Douglas \cite{PEERS} who introduced and analyzed a
low-order method. Later, Stenberg \cite{family} introduced a two and
three dimensional family. These method all use standard polynomial
basis function of Raviart-Thomas and Brezzi-Douglas-Marini type with
the addition of divergence free "bubble" degrees of freedom added in
order to enforce the stability of the rotational Lagrange multiplier. 

The families of the work \cite{family} where reconsidered by Cockburn,
Golapakrishnan and Guzm\' an \cite{CGG, GG} and they observed that
the stabilizing bubble spaces can be reduced. For the two dimensional
family this follows by inspecting which degrees of freedom  are
actually used in the proof given in \cite{family}. In three dimensions,
however, their contribution is significant. They introduced a tensor
valued bubble function, the use of which reduced the polynomial
degree of the divergence free degrees of freedom with one. 
 
 In \cite{family} and \cite{GG} the polynomials used for the stress
 and rotation are of the same degree with  additional bubbles of one
 degree higher for the stresses. However, if the degree of the
 rotation is lowered by one, then the bubbles needed are already
 contained in the space for the stress. This leads to the family of
 Arnold, Falk and Winther \cite{MR2249345,AFW07}. It however, leads to
 a decreased convergence rate for the stress. This lead us to
 introduce  a method in which the normal  component of the stress is
 reduced by one degree on each element edge or face. 
 
 In this paper we will perform a unified a priori error analysis
 covering the methods mentioned above. We essentially use the same
 analysis as in \cite{family}. In that paper we used the macroelement
 technique to establish the stability. The idea in this technique is
 that the crucial requirement for the stability is that the degrees of
 freedom are such that the local uniqueness (modulo a rigid body
 motion) is valid. The  local stability condition then follows from
 homogeneity  and continuity arguments. Below, we will, however,
 perform the analysis explicitly.
 
 In addition to the a priori analysis we will perform a novel a
 posteriori analysis. By modifying the classical Prager-Synge
 hypercircle theorem \cite{MR25902, NecasH}  we derive an estimate
 without any unknown constants. This estimate deteriorates for nearly
 incompressible materials and for that case we introduce another
 estimate.
 
 The plan of the paper is as follows. In the next section we recall
 the equations of linear elasticity with the rotational equilibrium as
 an independent equation. We prove the stability of the problem using
 the physical energy norms. In the next section we recall the finite
 element methods and prove the stability and a priori error estimates.
 In Section 4 we introduce and analyze a post processing of the
 displacement variable. This is crucial for the a posteriori analysis
 which is done in Section 5 and Section 6. In the last section we
 perform  numerical benchmark studies.

We use the established notation for Sobolev spaces and finite element
methods. We write  $A \lesssim B$ when there exists a positive
constant $C$,   that is independent of the mesh parameter and,
\emph{in particular}, of the two Lam\'e parameters $\mu, \lambda$ (see
below) such that $A \le C B$.  Analogously we define $A \gtrsim B$.
$A \lesssim B$ and $A \gtrsim B$ is denoted by $A \simeq B$. We are
careful of explicitly including the right physical parameters, and
this means that the dependency of the two Lam\'e parameters are made
explicit in the norms used.

\section{The equations of elasticity with rotational equilibrium}

Let $\Omega\subset \real^\dim$ be a polygonal or polyhedral domain.
The physical unknowns are the displacement vector
$u=(u_1,\dots,u_\dim)$ and the   stress tensor $\sigma=\{
\sigma_{ij}\}$, $i,j=1,\dots, \dim$. The loading consists of a body
load $f$ and a traction $g$ on the boundary part $\Gamma_N$. On the
complementary part $\Gamma_D$ homogeneous Dirichlet conditions for the
displacement are given.

The condition of force and rotational equilibrium gives the equations
\begin{equation}
\div \bfs+f=0 \  \mbox{ in  }  \Omega,
\end{equation}
and
\begin{equation}\label{roteq}
\bfs-\bfs^T=0 \  \mbox{ in  }  \Omega,
\end{equation}
respectively. 

The deformation is given by the infinitesimal strain tensor
\begin{equation}
\varepsilon(u) =\frac{1}{2} \big( \nabla u + \nabla u ^T\big).
\end{equation}The stress and strain tensors are related by a linear constitutive 
 law. For an isotropic material and the \emph{plane strain} ($\dim=2$)
 or \emph{three dimensional} ($\dim=3$)  problem this is given by
 \begin{equation}
 \comp \bfs =\varepsilon(u), 
 \end{equation}
 with the compliance tensor 
\begin{equation}\label{consteq1}
\comp \bft = \frac{1}{2\mu} \Big(  \bft -\frac{\lambda}{ 2\mu +\dim\lambda} \tr(\bft)I \Big),
\end{equation}
 where 
 $\mu \mbox{ and } \lambda$ are the  Lam\'e parameters.
 For the 
 \emph{plane stress} problem the  relationship is, cf.  \cite{NecasH},
 \begin{equation}\label{planestress}
  \comp \bft = \frac{1}{2\mu} \Big(  \bft -\frac{\lambda^*}{ 2\mu +2\lambda^*} \tr(\bft)I \Big),
\end{equation}
with
\begin{equation}
\lambda^* =\frac{2\mu \lambda}{\lambda +2\mu}.
\end{equation}
 
The inverse of the compliance matrix,  the elasticity matrix, we
 denote by $\elas=\comp^{-1}$, i.e. 
 \begin{equation} \label{eq::inccompl}
 \elas \tau =2\mu \tau +\lambda \tr(\tau) I,
 \end{equation} 
 for the plane strain and three dimensional problem, 
 and 
  \begin{equation} 
 \elas \bft = 2\mu \bft + \lambda^* \tr(\bft) I,
\end{equation}
for plane stress.
 
Due to the symmetry of the compliance tensor, the symmetry of the
stress tensor is implied by \eqref{consteq1}, \eqref{planestress},
and the traditional form of the elasticity equations are
\begin{eqnarray}
\comp \bfs -\varepsilon( u )&=&0\  \mbox{ in  }  \Omega,  \label{aaa}\\
\div \bfs + f&=&0 \  \mbox{ in  }  \Omega, 
\\
u&=&0 \  \mbox{ on } \Gamma_D,
\\
\bfs n&=&g \ \mbox{ on } \Gamma_N. \label{neumann} \end{eqnarray}
Using this formulation as a basis for mixed finite element methods leads, however, to rather complicated methods, cf. \cite{JM,ADG,AAW3D,MR3149075}, and this lead  Fraeijs de Veubeke to suggest mixed methods in which the rotational equilibrium \eqref{roteq} is treated as an independent equation \cite{FdV-W,FdV2}. Due to the introduction of one additional equation, one additional unknown is needed, and this is the skew symmetric infinitesimal rotation, viz.
\begin{equation}\label{rotdef} 
\rho = \frac{1}{2} \big( \nabla u -\nabla u^T\big).
\end{equation}
 This gives the set of 
 equations  
\begin{eqnarray}
\comp \bfs +\rho -\nabla u &=&0\  \mbox{ in  }  \Omega, \label{aaa-new}\\
  \bfs-  \bfs^T&=&0\  \mbox{ in  }  \Omega,\\
\div \bfs + f&=&0 \  \mbox{ in  }  \Omega, \label{e-eqs}
\\
u&=&0 \  \mbox{ on } \Gamma_D,
\\
\bfs n&=&g \ \mbox{ on } \Gamma_N. \label{neumann2} 
\end{eqnarray} 
The symmetric part of the first equation \eqref{aaa-new} yields the
traditional \eqref{aaa}, and the skew part gives \eqref{rotdef}.

In the sequel we will use the notation
\begin{equation}
\omega(u) = \frac{1}{2} \big( \nabla u -\nabla u^T\big).
\end{equation}
 The variational form of the problem is: find
$\bfs \in \HH$, $\rho\in \SK $ and $u\in \VV=\{\,  v\in [H^1(\Omega)]^\dim\, \vert \,
v\vert_{\Gamma_D}=0 \,\}$ such that 
\begin{equation}
\begin{aligned}
\B (\bfs,\rho,  u; \bft, \eta, v) &=(f,v) + \langle   g, v \rangle_{\Gamma_N} \\
& \forall(\bft,\eta, v) \in \HH \times \SK\times \VV,
\end{aligned}
\end{equation}
with the bilinear form 
\begin{equation}
\B (\bfs, \rho, u; \bft, \eta, v) =  (\comp \bfs,\bft )+ b( \bft; \rho, u)  -b(\bfs; \eta, v),
\end{equation}
where
\begin{equation}
b(\bft;\eta,v)= (\bft,\eta) -(\bft, \nabla v).
\end{equation}
The natural energy norms for analyzing this problem are 
 \begin{equation}
\enorm{ \bft } ^2= ( \comp \bft , \bft)\   \mbox{ and } \ \anorm{\varepsilon( v )}^2=   (\elas\varepsilon(v)  , \varepsilon(v)),
 \end{equation}
 which are the double of the strain energy expressed by the stress and displacement, respectively. The 
 Babu{\v{s}}ka--Brezzi condition is then simply the following identity.
 \begin{lemma} It holds that 
 \begin{equation}
 \begin{aligned}
\sup_{\bft\in \HH } \frac{\big(\bft, \varepsilon(v)\big)+(\eta-\omega(v), \bft)}{\enorm{\bft}} &= \big(\anorm{ \varepsilon(v)}^2+ \anorm{\eta-\omega(v)}^2\big)^{1/2}
\\&   \forall (\eta, v)\in \SK \times \VV.
\end{aligned}
\end{equation}
\end{lemma}
\begin{proof}For $ (\eta, v)\in \SK \times \VV$ given, we choose $\bft\in \HH$ by $\bft_{{\rm sym}} = {\mathcal A} \varepsilon(v) $ 
and $\bft_{{\rm skw}} ={\mathcal A}( \eta-\omega(v))$.  The claim then follows from the orthogonality between symmetric and skew symmetric tensors. 
\end{proof}

We denote
\begin{equation}
\E = \HH \times \SK \times \VV,
\end{equation}
and 
\begin{equation}\label{energynorm}
\norm{(\bft,\eta,v)}_\E^2 = \enorm{\bft}^2+ \anorm{\varepsilon(v) }^2 + \anorm{\eta- \omega(v) }^2.
\end{equation}

With the use of this energy norm, we first obtain the sharp "inf-sup" estimate above, and by using Lemma 3.1 of 
  \cite{HSV}, we then get the stability with a known constant.
  
\begin{theorem} \label{energystab} It holds that 
\begin{equation}
\begin{aligned}
 \sup_{(\varphi,\xi,z)
  \in \E}  \frac{\B(\bft,\eta,v; \varphi, \xi, z)} { \norm{(\varphi,\xi,z)}_\E
 } &\geq \Big(\frac{\sqrt{5}-1}{2} \Big)   \norm{(\bft,\eta,v)}_\E
 \\ & \qquad \qquad \forall (\bft,\eta,v)\in \E.
 \end{aligned}
 \end{equation}
 \end{theorem}
 For nearly incompressible materials the stability has to be posed differently, and we proceed as follows.
 For  $\bft$ we 
let $\bft^\D$ be the deviatoric part of $\bft$, defined by the condition $\tr(\bft^\D)=0. $ Hence, we have 
 \begin{equation}
  \bft = \bft^\D + \frac{1}{\dim} \tr(\bft) I.
  \end{equation}
  A direct computation gives.
\begin{lemma} \label{first} For the constitutive law \eqref{consteq1} it  holds that
  \begin{equation}\label{cbound}
  (\comp  \bft, \bft) = \frac{1}{2\mu} \Vert  \bft^\D\Vert_0^2 + \frac{1} {2\mu +\dim \lambda} \Vert \tr(\bft) \Vert_0^2.
  \end{equation}
  \end{lemma}
  From this we see that   $\enorm{ \cdot } $ ceases to define a norm in the incompressible limit $\lambda \to \infty$.
  In the general theory of Brezzi \cite{MR3097958}, the stability follows from the condition of the "ellipticity in the kernel". In \cite{LS1} we, however, gave an explicit proof.  In this part of the  stability estimate one has to use the 
  Babu{\v{s}}ka--Brezzi condition for the Stokes problem \cite{MR851383}.
    \begin{lemma} It holds that
    \begin{equation}\label{stokes}
  \sup_{v\in \VV} \frac{(\div v, q)}{\Vert  \varepsilon (v)\Vert_0} \geq \beta \Vert q \Vert_0 \quad \forall q \in L^2(\Omega).
  \end{equation}
    \end{lemma}  
The Babu{\v{s}}ka--Brezzi condition we now write as follows.
   \begin{lemma} \label{last} It holds that
    \begin{equation}
 \begin{aligned}
\sup_{\bft\in \HH } \frac{\big(\bft, \varepsilon(v)\big)+(\eta-\omega(v), \bft)}{\mu^{-1/2} \norm{\bft}_0}  &= \big(\mu  \norm{ \varepsilon(v)}_0^2+ \mu  \norm{\eta-\omega(v)}_0^2\big)^{1/2}\\&   \forall (\eta, v)\in \SK \times \VV.
\end{aligned}
\end{equation}
    \end{lemma}
    The full stability is a consequence of the above results. In the proof,  Lemma \ref{first} yields the stability for the deviatoric stress.  Inequality \eqref{stokes} implies the stability for the hydrostatic pressure part, whereas Lemma \ref{last} yields the control of the displacement and rotation. The proof is similar to that of Theorem 1 of \cite{LS1} and is omitted. 
    Since the stability constant will dependent on the unknown constant in \eqref{stokes}, there is no advantage to use an energy-type norm. By Korn's inequality it holds 
    \begin{align*}
\norm{ \varepsilon(v)}_0+ \norm{\eta-\omega(v)}_0 \simeq \norm{ \varepsilon(v)}+ \norm{\eta }_0,
\end{align*}
and using this we give the analysis with the norm
\begin{equation}\label{trinorm}
\Norm{(\bft,\eta,v)} ^2 = \mu^{-1}\norm{\bft}_0^2+ \mu \norm{\varepsilon(v) }_0^2 +\mu \norm{\eta}_0^2.
\end{equation}

   
\begin{theorem} \label{incompstab}It holds that 
\begin{equation}
\begin{aligned}
 \sup_{(\varphi,\xi,v)
  \in \E} \frac{\B(\bft,\eta,v; \varphi, \xi, z)} { \Norm{(\varphi,\xi,z)} 
 }& \gtrsim \Norm{(\bft,\eta,v)}
 \quad  \forall (\bft,\eta,v)\in \E.
 \end{aligned}
 \end{equation}
 \end{theorem}

 The finite element methods are based on a mixed formulation obtained by dualization. The space for the stress is 
\begin{equation}
\Hdiv= \{ \, \bft\in\HH\, \vert \, \div \bft \in [L^2(\Omega)]^\dim \, \}, 
\end{equation} that of 
the displacement is $\Vtwo$ and the rotation is in $\SK$.
The bilinear form used is  
 \begin{equation}
\M (\bfs, \rho , u; \bft, \eta, v) =  (\comp \bfs,\bft )+ b( \bft; \rho, u)  +b(\bfs; \eta, v),
\end{equation} 
with 
\begin{equation}\label{mixedb}
b( \bft; \rho, u) = (\div \bft, u)+(\bft,\rho).
\end{equation}
 The variational formulation is: find    $(\bfs, \rho, u)\in  \Hdg\times \SK \times [L^2(\Omega) ]^d$  such that
\begin{equation}\label{2ndvariational}
\begin{aligned}
\M (\bfs, \rho, u&; \bft,\eta,  v)+(f,v) =0
\\& \forall (\bft,\eta ,v) 
\in \Hdnull\times \SK \times [L^2(\Omega) ]^d,
\end{aligned}\end{equation}
with 
 \begin{equation}
 \Hdg= \{ \, \bft \in \H \, \vert \  \bft n \vert_{\Gamma_N} =g\, \},
\end{equation}
and
\begin{equation}
 \Hdnull= \{ \, \bft \in \H \, \vert \  \bft n \vert_{\Gamma_N} =0\, \} .
\end{equation}
 
\section{Weakly symmetric mixed finite element methods}
\label{sec::fem} 
 
The finite element methods are based on the variational formulation  \eqref{2ndvariational} posed in the subspaces $S_h\subset \Hdiv$, $R_h\subset \SK$ and $V_h\subset [L^2(\Omega) ]^d$:

\medskip

Find $\bfs_h \in S_h^g$, $\rho_h \in R_h$ and $u_h \in V_h$ such that
\begin{equation} \label{discsol}
\M (\bfs_h, \rho_h, u_h; \bft,\eta,  v)+(f,v) =0\quad \forall (\bft,\eta ,v) \in S_h^0 \times R_h \times V_h ,
\end{equation}
where 
\begin{equation}
S_h^g=\{ \bft \in S_h \, \vert \bft n =g_h  \, \mbox{ on } \Gamma_N\},
\end{equation}with $g_h$ being an approximation of $g$ (see below) and $S_h^0 =S_h \cap \Hdnull   $.

\medskip
 
We will give an analysis covering the methods mentioned in the introduction, i.e. those of Stenberg \cite{family}, Cockburn, Golapakrishnan and Guzm\' an \cite{CGG, GG}, and Arnold, Falk and Winther \cite{MR2249345,AFW07}. In addition, we will introduce a new method obtained by inspecting the one of Arnold et al. 

The finite element partitioning consists of  triangles or tetrahedrons and is denoted by $\Ch$. An edge (for $d=2$) or face ($d=3$) of an element $K\in \Ch$ we denote by $E$, and the collection of edges/faces is denoted by $\Gamma_h$.  By $P_l(K) $ and $\tilde P_l(K)$, with $l\geq0$,  we denote the polynomials and homogeneous polynomials of   degree $l$, respectively.

For a triangle the scalar valued bubble function $b_K$ is defined as
 \begin{equation}
 b_K = \lambda_1\lambda_2\lambda_3,
 \end{equation}
 where $\lambda_1, \, \lambda_2, \, \lambda_3$ are the barycentric coordinates. 
 For tetrahedrons the tensor valued bubble function is
 \begin{equation}
 b_K = \sum_{l=0}^3 \lambda_{l-3}\lambda_{l-2} \lambda_{l-1} \nabla \lambda_l^T\nabla  \lambda_l,
 \end{equation}
 where the sum is calculated modulo four,   and $\nabla  \lambda_l$ is considered as a row vector. In both dimensions it holds that the bubble is a third degree polynomial. The bubble induces a weighted inner product on $K\in \Ch$, viz.
 \begin{equation}
 (w b_K, v)_K = (v b_K, w)_K, \   (v b_K, v)_K\geq 0, \mbox{ and } (v b_K, v)_K=0 \ \Rightarrow v=0.
 \end{equation}
 With the aid of the bubble we define
\begin{equation}
   \Ml= \{ \, \bft   \, \vert\, \bft =\curl \big(\curl  (z)b_K), \,  z\in \Rl\, \},
\end{equation}
for $l=k-1$ and $l=k$. 
Here the $\curl$ is applied row wise to the tensors. The functions in
this space have both vanishing normal components at the boundary and a
vanishing divergence, viz.
\begin{equation}
\bft n \vert_{ \partial K} =0,  \quad \div \bft =0, \quad \forall \bft \in M_{l+1} (K) .
\end{equation} 
We  furter note that by the definition of the bubble function $b_K$ the
index of the space $\Ml$ refers to the polynomial order $l+1$ of the
corresponding functions for both $d=2$ and $d = 3$. The Raviart-Thomas
space on $K$ is
\begin{equation}
\RTk= [P_{k-1}(K)]^d \oplus x \tilde P_{k-1}(K).
\end{equation}
The tensors with the rows in $\RTk$ is denoted by
  \begin{equation}
[\RTk]^d .
\end{equation}

Next we   define the methods to be analyzed.

\begin{itemize}
\item

The \emph{finite element space for the displacement }is the same for all families, and we index it with $k\geq 1$:
\begin{equation}\label{vspace}
V_h=  \{\, v\in \Vtwo\, \vert \ v\vert _K\in [P_{k-1}(K)]^\dim\ \forall K\in \Ch\}.
\end{equation}

\item

In the  \emph{Cockburn-Golapakrishnan-Guzm\' an family } (CGG)
\cite{CGG}  the space for the rotation is defined by piecewise
polynomials of degree $k-1 \geq 1$
\begin{equation}
R_h =\{ \, \eta \in \SK \ \vert \ \eta\vert_K \in [P_{k-1}(K)]^\dd_{{\rm skw}} \ \forall K \in\Ch \},
\end{equation}
and that of the stress uses the Raviart-Thomas spaces and the bubble
space defined above
\begin{equation}
S_h = \{ \, \bft\in \Hdiv\, \vert \, \bft\vert_K \in [\RTk]^d + M_k(K) \  \forall K\in \Ch \}.
\end{equation}
In the error analysis below we will show that for a smooth solution this method yields the convergence rate
\begin{equation}
\Vert \bfs-\bfs_h\Vert_0 \lesssim h^k,
\end{equation}
which is optimal, i.e. the same as the interpolation with this approximation   space for the stress.
 \item

The \emph{Golapakrishnan and Guzm\' an} \cite{GG} refinement of
\emph{Stenbergs method} \cite{family} consists of the above choice
of $V_h$, and 
\begin{equation}
  R_h =\{ \, \eta \in \SK \ \vert \ \eta\vert_K \in [P_{k}(K)]^\dd_{{\rm skw}} \ \forall K \in\Ch \},
  \end{equation}
  and 
\begin{equation}\label{sspace_sgg}
S_h= \{ \, \bft\in \Hdiv\, \vert \, \bft\vert_K \in [P_k(K)]^{d \times d} + \Mk \  \forall K \in \Ch\},
\end{equation}
with $k\geq 1.$
Below the family is referred to as SGG. 
The convergence rate is now also optimal
\begin{equation}
\Vert \bfs-\bfs_h\Vert_0 \lesssim h^{k+1}.
\end{equation}
\item

In the \emph{Arnold-Falk-Winther family} (AFW) \cite{MR2249345,AFW07} the polynomial degree for the rotation is lowered by one, and a consequence of that, is that the additional bubble degrees of freedom are not needed. The spaces are then defined for $k\ge 1$ by
 \begin{equation}\label{sspace_afw}
S_h= \{ \, \bft\in \Hdiv\, \vert \, \bft\vert_K \in [P_k(K)]^{d \times d} \  \forall K\in \Ch \},
\end{equation}
and
\begin{equation}\label{rspace_afw} 
R_h =\{ \, \eta \in \SK \ \vert \ \eta\vert_K \in  [P_{k-1}(K)]^\dd_{{\rm skw}}\ \forall K \in\Ch \}.
\end{equation}
The drawback of this choice is that the polynomials chosen for the
 unknowns are not in balance, and as a consequence the convergence rate is not optimal. It holds 
 \begin{equation}\label{sameacc}
\Vert \bfs-\bfs_h\Vert_0 \lesssim h^{k},
\end{equation}
which should be compared with
 \begin{equation}
\inf_{\bft\in S_h} \Vert \bfs-\bft\Vert_0 \lesssim h^{k+1}.
\end{equation}
 
  \item
 This motivates the following modification which appears to be new. On each inter element boundary the polynomial order is reduced by one. We set
  \begin{equation}
 S_k^r(K) = \{ \, v\in [P_k(K)]^\dd\, \vert \ vn\vert_ E  \in [P_{k-1} (E)]^d \ \forall E\subset \partial K \}
 \end{equation}
and define the stress space  for $k\geq 2$ as
 \begin{equation}\label{sspace_rafw}
S_h= \{ \, \bft\in \Hdiv\, \vert \, \bft\vert_K \in  S_k^r(K) \  \forall K\in \Ch \}.
\end{equation}
The space for the rotation is as in the original method. This method
we will refer to as $\NEW$ (reduced AFW).

As the most
 efficient method of solving the finite element system is by
 hybridization \cite{AB}, the decrease of the normal degrees of
 freedom on element boundaries leads to a decrease of computational
 cost compared with the corresponding AFW method, but with the same accuracy, i.e., \eqref{sameacc}. Most drastically this is seen for the lowest order method, with $k=2$, for which the numbers of unknowns are halved.
\end{itemize}

 In the methods  $g_h\vert_E =Q_E g\,$ with the $L^2$ projection $Q_E:
 [ L^2(E)]^d \to [P_{l} (E)]^d$, with $l=k$ for  SGG and AFW, and
 $l=k-1$ for the others. Further,  we define $Q_h$ by $Q_h|_E = Q_E$.

In the error analysis we will use the broken energy norm for the displacement 
\begin{equation}
\hnorm{v}= \big(  \sum_{K\in \Ch}   \Vert \varepsilon(v) \Vert_{0,K} ^2+ \sum_{E\in \Gamma_{h }}h_E ^{-1} 
\Vert
\jump{v}\Vert_{0,E}^2+ \sum_{E\subset \Gamma_D }h_E ^{-1} 
\Vert
 v\Vert_{0,E}^2
  \big)^{1/2}.
\end{equation}
In the proof of the discrete   Babu{\v{s}}ka--Brezzi  stability condition we first obtain the stability in the norm
\begin{equation}
\begin{aligned}
  \hnorm{v}^2+ \sum_{K\in \Ch} \norm{\eta-\omega(v) }_{0,K}^2\quad  \forall (v,\eta) \in V_h\times R_h .
 \end{aligned}
\end{equation}

From the  discrete Korn inequality \cite{brennerkorn} it follows that
\begin{equation} \hnorm{v}^2 \gtrsim \sum_{K\in \Ch} \Vert \omega(v) \Vert_{0,K}^2  \quad \forall  v\in V_h,
\end{equation}

This implies the norm equivalence.
\begin{lemma}\label{normeq} It holds that
\begin{equation}
  \hnorm{v}^2+ \sum_{K\in \Ch} \norm{\eta-\omega(v) }_{0,K}^2\simeq  \hnorm{v}^2+ \norm{\eta }_{0}^2  \quad  \forall (v,\eta) \in V_h\times R_h .
 \end{equation}
\end{lemma}
Hence, we will use the norm
\begin{equation}
\Normh{(\bft, \eta,v) }^2 =  \mu^{-1}\norm{\bft}_0^2+ \mu \hnorm{v}^2 +\mu \norm{\eta}_0^2.
\end{equation}
The discrete  Babu\v{s}ka-Brezzi is the following.
    \begin{lemma} \label{discinfsup} Let  $b$ be the bilinear form \eqref{mixedb}. It holds that
\begin{equation} \label{stability2}
\sup_{\bft\in\Sh^0}  \frac{b( \bft; \eta, v) }
{ \norm{\bft}_0}
 \gtrsim  (\norm{\eta } _0 +\norm{v}_\etth)  \quad \forall (v,\eta) \in \Vh\times R_h.
\end{equation}
\end{lemma}
In the proof we proceed differently for the cases $k\geq 2$ and $k=1$. For the first case we need an additional result.
  Define 
 \begin{equation}
 R_h^0 = \{ \, \eta \in R_h\, \vert \, \eta\vert_K \in [P_0(K)]^\dd_{{\rm skw}} \ \forall K\in \Ch\, \},
 \end{equation}
 and let $N_h : R_h \to R_h^0$ be the $L^2$ projection.
\begin{lemma}\label{bubbstab}
It holds that
\begin{equation}
\sup_{\bft \in \Ml  } \frac{(\bft, \eta)_K}{\norm{\bft}_{0,K}} \gtrsim \norm{\eta-\Nh \eta}_{0,K} \quad \forall \eta \in \Rl.
\end{equation}
\end{lemma}
\begin{proof} For $\eta \in  \Rl$ given, we choose 
  $\bft\in \Ml  $ by
\begin{equation}
\bft \vert_K = h_K^2 \curl \big(\curl  (\eta)b_K\big) .
\end{equation}
Integration by parts gives
\begin{equation}
(\bft, \eta)_K = h_K ^2( \curl (\eta )b_K, \curl \eta)_K \gtrsim h_K ^2\Vert \curl \eta\Vert_{0,K}^2. 
\end{equation}
Since, $\eta $ is skew symmetric
it holds
\begin{equation}
\Vert \eta-\Nh\eta  \Vert_{0,K} \lesssim h_K \Vert \nabla  \eta\Vert_{0,K}\lesssim h_K \Vert \curl \eta\Vert_{0,K}.
\end{equation}
By the inverse inequality we have 
\begin{equation}
\begin{aligned}
\Vert \bft  \Vert_{0,K}^2 &= h_K^2 \Vert \curl \big(\curl  (\eta)b_K\big)\Vert_{0,K} ^2=
h_K^2 \Vert \curl \big(\curl  (\eta-\Nh\eta )b_K\big)\Vert_{0,K} ^2
\\
&
 \lesssim \Vert \eta-\Nh\eta  \Vert_{0,K}^2  .
\end{aligned}
\end{equation}
 \end{proof}

Using  this result we give the 

\noindent {\em Proof of Lemma \ref{discinfsup}} for $k\geq 2$. 
We first note that $V_h$ is the same for all methods considered. 
We define
\begin{equation}\label{RTdofs}
\RTkH =\{ \, \bft\in \RTk\, \vert \, \, \bft \cdot n\vert_E \in P_1(E) \  \forall E \subset \partial K\,   \}.
\end{equation}
The degrees of freedom of this space are \cite{MR3097958} 
\begin{equation}
\label{mindofs}
\int_K \bft \cdot v, \ v\in  [P_{k-2}(K)]^d   \ \mbox{ and } \int_E( \bft \cdot n) v, \ v\in P_1(E), \ E \subset \partial K.
\end{equation}
By $Q_E^1:[ L^2(E) ]^{\dim} :\to [ P_1(E) ]^{\dim} $ we denote the $L^2$ projection. 

Let $(v,\eta) \in V_h\times R_h$ be given. For $\bft_1 \in S_h$, with $\bft_1\vert_K \in [\RTkH]^d\subset S_h\vert_K $, it holds
\begin{equation}
 \begin{aligned}
b(\bft_1;& v, \eta)  = (\div \bft_1, v)+(\bft_1,\eta)
\\
=  &- \sum_{K\in \Ch} (\bft_1, \nabla v)_K 
+\sum_{E\in \Gamma_h } \int_E( \bft_1n) \! \cdot\! \jump{v} +\sum_{E\subset \Gamma_D } \int_E( \bft_1n)\! \cdot \! v +
( \bft_1,\eta)
\\
=&- \sum_{K\in \Ch} (\bft_1, \nabla v)_K +\sum_{E\in \Gamma_h } \int_E( \bft_1n) \! \cdot \!\jump{Q_E^1v} +\sum_{E\subset \Gamma_D } \int_E( \bft_1n) \! \cdot \! Q_E^1v 
+ ( \bft_1,\eta)
\\
=& \sum_{K\in \Ch} (\bft_1,\Nh\eta- \nabla v)_K +\sum_{E\in \Gamma_h } \int_E( \bft_1n) \! \cdot \! \jump{Q_E^1v} +\sum_{E\subset \Gamma_D } \int_E( \bft_1n) \! \cdot \! Q_E^1v 
 \\
&
+(\bft_1, \eta-\Nh \eta).
\end{aligned}
\end{equation}
By \eqref{mindofs} we can now choose $\bft_1$ such that 
\begin{equation}\label{xx}
\begin{aligned}
  &\sum_{K\in \Ch} (\bft_1,\Nh \eta- \nabla v)_K+\sum_{E\in \Gamma_h } \int_E( \bft_1n) \! \cdot \! \jump{v} +\sum_{E\subset \Gamma_D } \int_E( \bft_1n) \! \cdot \! v
  \\
= &\sum_{K\in \Ch} \Vert \Nh \eta- \nabla v\Vert_{0,K}^2+\sum_{E\in \Gamma_h }h_E^{-1} \int_E \vert \jump{Q^1_Ev} \vert^2+\sum_{E\subset \Gamma_D} h_E^{-1}\int_E\vert  Q^1_E v\vert^2,
\end{aligned}
\end{equation}
and  
\begin{equation}\label{xy}
\Vert \bft_1\Vert_0 \leq C_1 \Big(  \sum_{K\in \Ch} \Vert \Nh \eta- \nabla v\Vert_{0,K}^2+\sum_{E\in \Gamma_h }h_E^{-1} \int_E \vert \jump{Q^1_Ev} \vert^2+\sum_{E\subset \Gamma_N }h_E^{-1} \int_E\vert  Q^1_E v\vert^2      \Big)^{1/2}.
\end{equation}
Note that $\Ml \subset S_h\vert_K$, with $l=k-1$ for CGG, AFW, \NEW, and $l=k$ for SGG. 
Using Lemma \ref{bubbstab} we then choose $ \bft_2|_K \in \Ml$,  such that  $\div
\bft_2=0$,  
\begin{equation}\label{yy}
(\bft_2, \eta)   = \Vert   \eta-\Nh  \eta\Vert_{0}^2,  \quad \mbox{ and } \quad 
\Vert \bft_2 \Vert_{0}\leq C_2 \Vert \eta-\Nh\eta  \Vert_{0}.
\end{equation}
 For $\bft =\delta \bft_1+\bft_2$, with $\delta>0$, it then holds
\begin{equation}
 \begin{aligned}
b(\bft; v, \eta)  =& (\div \bft, v)+(\bft,\eta)= (\div \delta \bft_1, v)+(\delta \bft_1+\bft_2,\eta)
\\
= &\delta\Big(- \sum_{K\in \Ch} (\bft_1, \nabla v)_K +\sum_{E\in \Gamma_h } \int_E( \bft_1n) \! \cdot \! \jump{v} +\sum_{E\subset \Gamma_D } \int_E( \bft_1n) \! \cdot \! v\Big)
\\& +
(\delta \bft_1+\bft_2,\eta)
\\ 
=&\delta\Big(
\sum_{K\in \Ch} (\bft_1,N_h\eta- \nabla v)_K +\sum_{E\in \Gamma_h } \int_E( \bft_1n) \! \cdot \! \jump{v} +\sum_{E\subset \Gamma_D } \int_E( \bft_1n) \! \cdot \! v\Big) 
\\&+ \delta (\bft_1,\eta -\Nh \eta) +
(\bft_2,\eta).
\end{aligned}
\end{equation}
 By  Schwarz and the arithmetic geometric mean inequality we have
 \begin{equation}
  \delta (\bft_1,\eta -\Nh \eta)\geq - \delta \norm{\bft_1}_0 \norm{\eta -\Nh \eta}_0 \geq-  \frac{\delta^2} {2}  \norm{\bft_1}_0^2 -\frac{1}{2}  \norm{\eta -\Nh \eta}_0^2.
 \end{equation}
 Hence, using \eqref{xx}, \eqref{xy}, and \eqref{yy},  it holds 
 \begin{equation}
 \begin{aligned}
&b(\bft;   v, \eta) \geq
\\
  &\delta\big(1-\frac{\delta C_1}{2}\big) \Big(  \sum_{K\in \Ch} \Vert \Nh \eta- \nabla v\Vert_{0,K}^2+\sum_{E\in \Gamma_h } h_E^{-1}\int_E \vert \jump{Q^1_Ev} \vert^2+\sum_{E\subset \Gamma_D }h_E^{-1} \int_E\vert  Q^1_E v\vert^2\Big)    
\\& +\frac{1}{2}  \Vert   \eta-\Nh  \eta\Vert_{0}^2. 
 \end{aligned}
 \end{equation}
When choosing $0<\delta < 2C_1^{-1} $, we have
\begin{equation}
 \begin{aligned}
b(\bft;  v, \eta) \gtrsim    \sum_{K\in \Ch} &\Vert \Nh \eta- \nabla v\Vert_{0,K}^2+\sum_{E\in \Gamma_h }h_E^{-1} \int_E \vert \jump{Q^1_Ev} \vert^2
\\
& +\sum_{E\subset \Gamma_D } h_E^{-1}\int_E\vert  Q^1_E v\vert^2\ + \Vert   \eta-\Nh  \eta\Vert_{0}^2,
 \end{aligned}
 \end{equation}
 and
 \begin{equation}
  \begin{aligned}
 \norm{\bft }_0 \lesssim    \Big( \sum_{K\in \Ch} &\Vert \Nh \eta- \nabla v\Vert_{0,K}^2+\sum_{E\in \Gamma_h }h_E^{-1} \int_E \vert \jump{Q^1_Ev} \vert^2
 \\& +\sum_{E\subset \Gamma_D } h_E^{-1}\int_E\vert  Q^1_E v\vert^2     
  + \Vert   \eta-\Nh  \eta\Vert_{0}^2\Big)^{1/2}. 
\end{aligned}
 \end{equation}
  By the orthogonality between symmetric and skew symmetric tensors we get
 \begin{equation}
 \Vert \Nh \eta- \nabla v\Vert_{0,K}^2= \Vert \Nh \eta- \omega(v)\Vert_{0,K}^2+\Vert   \varepsilon( v)\Vert_{0,K}^2.
 \end{equation}
The triangle inequality yields
\begin{equation}
\begin{aligned}
&\Vert   \eta-\Nh  \eta\Vert_{0,K}^2 +\Vert \Nh \eta- \omega(v)\Vert_{0,K}^2+\Vert   \varepsilon( v)\Vert_{0,K}^2
\\&\gtrsim \Vert   \eta- \omega(v)\Vert_{0,K}^2+\Vert   \varepsilon( v)\Vert_{0,K}^2,
\end{aligned}
\end{equation}
and by scaling we have
\begin{equation}
\begin{aligned}
&\sum_{K\in \Ch} \Vert   \varepsilon( v)\Vert_{0,K}^2+\sum_{E\in \Gamma_h } h_E^{-1}\int_E \vert \jump{Q^1_Ev} \vert^2+\sum_{E\subset \Gamma_D } h_E^{-1}\int_E\vert  Q^1_E v\vert^2 
\\
& \simeq 
\sum_{K\in \Ch} \Vert   \varepsilon( v)\Vert_{0,K}^2+\sum_{E\in \Gamma_h }h_E^{-1} \int_E \vert \jump{ v} \vert^2+\sum_{E\subset \Gamma_D } h_E^{-1}\int_E\vert   v\vert^2.
\end{aligned}
\end{equation}
Combining the above estimates proves that
\begin{equation}  
\sup_{\bft\in\Sh^0}  \frac{b( \bft; \eta, v) }
{ \norm{\bft}_0}
 \gtrsim \big( \  \hnorm{v}^2+ \sum_{K\in \Ch} \norm{\eta-\omega(v) }_{0,K}^2\big)^{17/2} \quad \forall (v,\eta) \in \Vh\times R_h,
\end{equation}
and the 
assertion follows from the norm equivalence of Lemma \ref{normeq}.
 
$\square$

Next, we consider the low order methods, the families SGG and AFW with
$k=1$, i.e. methods with piecewise constants for the displacement.
The rigid body motions on $K$ we denote by
   \begin{equation}
   RM(K) = \begin{cases}  
   \{ \, v \, \vert \, v=(a,b) +c(-x_2, x_1) , \ a,b,c \in \real \, \} &\mbox{ for } d=2,
   \\   \{ \, v \, \vert \, v=a + b\times x, \, a,b \in \real^3 \, \} &\mbox{ for } d=3.
    \end{cases}
   \end{equation}
 With this we define
   \begin{equation}
   RM_h= \{ \, \psi\,  \vert \, \psi\vert_K \in RM(K)\  \forall K\in \Ch\, \}.
   \end{equation}
   
   For $ \Nh\eta\in R_h^0$ (note that $\eta$ is linear for SGG with $k=1$) there
   exist a unique $r( \Nh\eta)\in RM_h $  such that  
 \begin{align}
 N_h   \eta\vert_K &= \omega\big(r( \Nh\eta\vert_K )\big)= \nabla r( \Nh\eta\vert_K),  \label{help1}
\\ 
 \mbox{ and } & \int_K r( \Nh\eta)=0 \quad \forall K\in \Ch. \label{help2}
 \end{align}
 
 Let us now give the
 
 \noindent {\em Proof of Lemma \ref{discinfsup} for the families SGG and AFW, with $k=1$.}
 
 Consider first AFW with 
  \begin{equation}\label{sspace}
S_h= \{ \, \bft\in \Hdiv\, \vert \, \bft\vert_K \in [P_1(K)]^{d \times d} \  K\in \Ch\, \},
\end{equation}\ 
and
\begin{equation}\label{rspace} 
R_h =\{ \, \eta \in \SK \ \vert \ \eta\vert_K \in  [P_0(K)]^\dd_{{\rm skw}} \ \forall K \in\Ch\, \}.
\end{equation}
 Let $(v,\eta) \in V_h \times R_h$ be given. Now, $R_h^0=R_h$ and we write $\eta\vert  _K= \nabla (r(\eta|_K))$. 
  Since $\div \bft \vert_K \in P_0(K)$, \eqref{help2} gives
  \begin{equation}
 \begin{aligned}
 (\bft, \eta) &= \sum_{K\in \Ch} \big(\bft, \nabla r(\eta )\big)_K 
 \\& = -  \sum_{K\in \Ch} \big(\div \bft,  r(\eta )\big)_K  
  +\sum_{E\in \Gamma_h } \int_E( \bft n) \! \cdot \! \jump{ r(\eta )} +\sum_{E\subset \Gamma_D } \int_E( \bft n) \! \cdot \!r(\eta )
 \\& =  \sum_{E\in \Gamma_h } \int_E( \bft n)\! \cdot \! \jump{ r(\eta )} +\sum_{E\subset \Gamma_D } \int_E( \bft n) \! \cdot \! r(\eta ).
 \end{aligned}
  \end{equation}
On the other hand, for $(\div \bft, v)$ we have
   \begin{equation}
 \begin{aligned}
 (\div \bft, v)  &= \sum_{K\in \Ch}  (\div \bft, v)_K  
 \\& = -\sum_{K\in \Ch}  (  \bft, \nabla v)_K  + \sum_{E\in \Gamma_h } \int_E( \bft n) \! \cdot \! \jump{ v} +\sum_{E\subset \Gamma_D } \int_E( \bft n) \! \cdot \! v
  \\& =  \sum_{E\in \Gamma_h } \int_E( \bft n) \! \cdot \! \jump{ v} +\sum_{E\subset \Gamma_D } \int_E( \bft n) \! \cdot \! v.
 \end{aligned}
  \end{equation}
  Hence, 
    \begin{equation}
 \begin{aligned}
 b(\bft; \eta, v)&=(\div \bft, v) +(\bft, \eta)
 \\&
  =   \sum_{E\in \Gamma_h } \int_E( \bft n) \! \cdot \! \big(\jump{ r(\eta )+v} \big)+\sum_{E\subset \Gamma_D} \int_E( \bft n) \! \cdot \! \big(r(\eta )+v\big).
 \end{aligned}
  \end{equation}
  We now choose $\bft$  such that
  \begin{equation}
   \bft n= h_E^{-1} \big(\jump{ r(\eta )+v} )\ \mbox{ on } E\subset \Gamma_h, \ \mbox{ and } 
    \bft n= h_E^{-1} \big( r(\eta )+v )\ \mbox{ on } E\subset \Gamma_D,
  \end{equation}
  which gives 
  \begin{equation}
  \begin{aligned}
   b(\bft; \eta, v) =  
      \sum_{E\in \Gamma_h } h_E^{-1} \int_E \big(\jump{ r(\eta )+v} \big)^2+\sum_{E\subset \Gamma_D }h_E^{-1} \int_E \big(r(\eta )+v\big)^2
      ,
      \end{aligned}
  \end{equation}
  and by scaling
  \begin{equation}
 \norm{\bft}_0 \lesssim   \Big(\sum_{E\in \Gamma_h } h_E^{-1} \int_E \big(\jump{ r(\eta )+v} \big)^2+\sum_{E\subset \Gamma_D }h_E^{-1} \int_E \big(r(\eta )+v\big)^2\Big)^{1/2}.
  \end{equation}
 Since $ \varepsilon(r(\eta )+v)\vert_K =0$, the discrete Korn inequality \eqref{korn} gives
  \begin{equation}
  \norm{r(\eta )+v}_ \etth \lesssim  \Big(\sum_{E\in \Gamma_h } h_E^{-1} \int_E \big(\jump{ r(\eta )+v} \big)^2+\sum_{E\subset \Gamma_D }h_E^{-1} \int_E \big(r(\eta )+v\big)^2\Big)^{1/2}.
  \end{equation}
Next, \eqref{help1} gives
  \begin{equation}
  \begin{aligned}
  \norm{&r(\eta )+v}_ \etth^2 = \norm{\eta}_0^2 
  \\ &+  \Big(\sum_{E\in \Gamma_h } h_E^{-1} \int_E \big(\jump{ r(\eta )+v} \big)^2+\sum_{E\subset \Gamma_D }h_E^{-1} \int_E \big(r(\eta )+v\big)^2\Big).
  \end{aligned}
  \end{equation}
  By scaling arguments we get 
   \begin{equation}
  \begin{aligned}
 \norm{\eta}_0^2 
&+  \Big(\sum_{E\in \Gamma_h } h_E^{-1} \int_E \big(\jump{ r(\eta )+v} \big)^2+\sum_{E\subset \Gamma_D }h_E^{-1} \int_E \big(r(\eta )+v\big)^2\Big)\\
&\simeq  \norm{\eta}_0^2 +\norm{v}_ \etth^2,
  \end{aligned}
  \end{equation}
  and the assertion is established for AFW.
  
  For the lowest order SGG the finite element space for the deflection
  is the same, and that of the rotation is of one degree higher. In
  this case $\eta \in R_h$ is split as $\eta=\Nh \eta +(I-\Nh) \eta$.
  The stability of the part $ \norm{\Nh \eta}_0 +\norm{v}_ \etth$ is
  established as above, and for the part $\norm{(I-\Nh)\eta}_0$ the
  degrees of freedoms in $\Mtwo$ are used, following the same
  arguments as in the proof of Lemma \ref{discinfsup} for $k > 1$.
  $\square$

 In addition to the stability estimate \eqref{stability2}, the following   property is needed in order to get an optimal a priori error estimate
 \begin{equation}\label{eqprop}
 \div \bft \in V_h \quad \forall \bft \in S_h.
 \end{equation}
 From this it follows that
 \begin{equation}\label{projpro}
 (\div  \bft, v-P_hv) =0  \quad    \forall \bft\in S_h,
 \end{equation}
 where $P_h: [L^2(\Omega) ]^\dim \to V_h$
 is the $L^2$ projection.

To obtain a stability estimate uniformly valid with respect to the incompressibility we will use the following estimate.

\begin{lemma}\label{incuse} It holds that
\begin{equation}
\sup_{v\in  V_h} \frac{(v, \div \bft)}{\Vert v \Vert_\etth} \geq C_1\Vert \tr(\bft) \Vert_0-C_2\Vert \bft^\D \Vert_0 \quad \forall \bft\in \Sh^0.
\end{equation}
\end{lemma}

\begin{proof}
Given $\bft\in \Sh^0$, \eqref{stokes}  implies that there exists $z\in \VV$ such that
\begin{equation}
(\div v, \tr(\bft)) = -\beta \Vert \tr(\bft)\Vert_0^2 \ \mbox{and } \ \Vert \varepsilon(v) \Vert_0=  \Vert \tr(\bft)\Vert_0.
\end{equation}
Let $P_h v\in V_h$ be the projection in \eqref{projpro}.
It holds
\begin{equation}
\begin{aligned}
 (P_h v&, \div \bft) = (v, \div \bft)  
=-(  \varepsilon(v),  \bft) 
\\& =   - (  \div v,  \tr(\bft) )-(  \varepsilon(v),  \bft^d) 
   \geq \beta \Vert \tr(\bft)\Vert_0^2  -\Vert   \varepsilon(v) \Vert \Vert\bft^d\Vert_0
  \\& = \Vert \varepsilon(v) \Vert_0\Vert_0\big(  \beta \Vert \tr(\bft)\Vert_0    -   \Vert\bft^d\Vert_0\big).
\end{aligned}
\end{equation}
By scaling we have
\begin{equation}
\hnorm{P_h v} \lesssim \Vert \varepsilon(v) \Vert_0.
\end{equation}
Combining the two estimates above proves the claim.
\end{proof}

From Lemma \ref{first},  Lemma \ref{discinfsup}, and  Lemma \ref{incuse},   we get the stability of the method.

\begin{theorem} \label{error1} It holds that
 \begin{equation}\label{stab2}
\begin{aligned}
 \sup_{(\varphi,\xi,v)
  \in \Sh^0\times R_h \times\Vh} \frac{\M(\bft,\eta,v; \varphi, \xi, z)} { \Normh{(\varphi, \xi, z)}
 }& \gtrsim   \Normh{(\bft, \eta, v)}
 \\ & \qquad  \forall (\bft,\eta,v)\in \Sh^0\times R_h \times\Vh.
 \end{aligned}
 \end{equation}
 \end{theorem} 
  
In the proof of the error estimate we will use an averaging interpolation operator applied to the discrete displacement space.   For $k\geq 2$ it is defined as the so called Oswald interpolation $I_h^a :V_h\to V_h \cap [H^1_D(\Omega)]^\dim$ (cf. \cite[Chapter 5.5.2]{DiPietroErn}), which satisfies
\begin{equation}\label{oswald}
 \Vert   \nabla  \Ioh   v\Vert _0\ + \Big(  \sum_{K\in \Ch}   h_K^{-2} \Vert  v-\Ioh v\Vert _{0,K}^2 \Big)^{1/2} \lesssim \norm{v}_ \etth.
\end{equation}
For $k=1$ the displacement space is
\begin{equation} 
V_h^0=  \{\, v\in \Vtwo\, \vert \ v\vert _K\in [P_{0}(K)]^\dim\ \forall K\in \Ch\},
\end{equation}
and for this we define the operator $\Ioh:V_h^0 :\to V_h^1\cap [H_D(\Omega)]^\dim$, with
\begin{equation} 
V_h^1=  \{\, v\in \Vtwo\, \vert \ v\vert _K\in [P_{1}(K)]^\dim\ \forall K\in \Ch\},
\end{equation}
in the following way. The vertex values of $\Ioh v$ are the average of
the values of $v$ in the elements sharing that vertex, and letting
$\Ioh v=0$ on $\Gamma_D$. It is easily proved that this interpolation
satisfies \eqref{oswald}.

We will next derive the quasi-optimality of the methods.  Here $f_h$ is any piecewise polynomial approximation of $f$.
\begin{theorem} \label{Theo3} It holds that
\begin{equation}
\begin{aligned}
\mu^{-1/2} \norm{&\bfs-\bfs_h}_0 + \mu^{1/2}( \norm{P_h u-\ph}_ \etth + \norm{ \rho-\rho_h}_0) 
\\  \lesssim\,
\mu^{-1/2}& \big(\inf_{\tau \in \Sh^g }\Vert \bfs-\tau \Vert_0 +\big(  \sum_{K\in \Ch} h_K^2 \Vert f- f_h\Vert_{0,K}^2     \big)^{1/2}\big)
\\&+\inf_{\eta \in R_h}\mu^{1/2} \Vert \rho-\eta \Vert_0 .
  \end{aligned}
  \end{equation}
\end{theorem}
\begin{proof}  
  By Theorem  \ref{error1} there exist
  $(\varphi, \xi, v) \in \Sh^0\times R_h \times V_h$ with
  \begin{equation}\label{normal}
  \Normh{ (\varphi,\xi, v)} =1,
  \end{equation}
   such that for all
  $\bft \in \Sh^g$ and $\eta\in R_h$ we have
\begin{equation}
\begin{aligned}
\snorm{  \bfs_h-\bft}&+\mu^{1/2}(\norm{\rho_h-\eta}_0 + \norm{P_h u-u_h }_ \etth )
\\
 &
\lesssim \M(   \bfs_h-\bft, \rho_h-\eta, u_h -P_h u   ; \varphi, \xi  , v).
\end{aligned}
\end{equation}
 By the consistency we have
\begin{equation}
 \M(   \bfs_h-\bft, \rho_h-\eta, u_h -P_h u   ; \varphi, \xi  , v)=  \M(   \bfs-\bft, \rho-\eta, u -P_h u   ; \varphi, \xi  , v).
\end{equation}
Writing out gives 
\begin{equation}
\begin{aligned}
 \M(   \bfs-\bft,& \rho-\eta, u -P_h u    ; \varphi, \xi  , v)
 \\=  &(\comp( \bfs-\bft) ,  \varphi ) +  (\rho-\eta,\varphi  ) +(u -P_h u,\div \varphi) 
\\& +(\bfs-\bft, \xi)+  (\div (\bfs-\bft), v) .
\end{aligned}
\end{equation}
From \eqref{cbound} and \eqref{normal} it follows
\begin{equation}
 (\comp( \bfs-\bft),  \varphi)\lesssim \snorm{\bfs-\bft }  \,\snorm{   \varphi} \lesssim  \snorm{  \bfs-\bft  } ,
\end{equation}
\begin{equation}
 (\rho-\eta,\varphi  ) \leq\norm{\rho-\eta}_0\norm{\varphi }_0 \leq \mu^{1/2} \norm{\rho-\eta}_0\,\mu^{-1/2}  \norm{\varphi }_0  \leq \mu^{1/2} \norm{\rho-\eta}_0 ,  
\end{equation}
and
\begin{equation}
 (\bfs-\bft, \xi) \leq\norm{\bfs-\bft}_0\norm{\xi}_0 \leq\mu^{-1/2} \norm{\bfs-\bft}_0\, \mu^{1/2}\norm{\xi}_0 \leq\mu^{-1/2} \norm{\bfs-\bft}_0 .
\end{equation}
 By the property \eqref{eqprop} the second term vanishes
 \begin{equation}
  (u -P_h u,\div \eta) =0.
\end{equation}
Next, we write
\begin{equation}
(\div(\bfs-\bft), v) = (\div(\bfs-\bft),  v-\Ioh v)  +  (\div(\bfs-\bft), \Ioh v).
\end{equation}
The first term above we treat as follows. The relations \eqref{oswald} and \eqref{normal} yield
\begin{equation}
\begin{aligned}
 (\div&(\bfs-\bft),  v-\Ioh v)    =\sum_{K\in \Ch}  (\div(\bfs-\bft),  v-\Ioh v)_K
 \\
 & \leq 
 \sum_{K\in \Ch}  \Vert \div(\bfs-\bft)\Vert_{0,K} \Vert  v-\Ioh v\Vert _{0,K}
 \\ &\leq \Big( \mu^{-1} \sum_{K\in \Ch}  h_K^2 \Vert \div(\bfs-\bft)\Vert_{0,K}^2  \Big) ^{1/2} \Big(  \mu \sum_{K\in \Ch}   h_K^{-2} \Vert  v-\Ioh v\Vert _{0,K}^2 \Big)^{1/2}
 \\
& \lesssim
 \Big( \mu^{-1} \sum_{K\in \Ch}  h_K^2 \Vert \div(\bfs-\bft)\Vert_{0,K}^2  \Big) ^{1/2}.
 \end{aligned}
\end{equation}
By a posteriori error analysis techniques \cite{MR3059294} we have
\begin{equation}
h_K \Vert \div(\bfs-\bft)\Vert_{0,K} \lesssim \big(\Vert \bfs-\bft\Vert_{0,K}+ h_K \Vert f-f_h\Vert_{0,K}\big)  ,
\end{equation}
and hence
\begin{equation}
 (\div(\bfs-\bft),  v-\Ioh v) \lesssim
  \mu^{-1/2} \big(\Vert  \bfs-\tau \Vert_0  +\big(  \sum_{K\in \Ch} h_K^2 \Vert f- f_h\Vert_{0,K}^2\big)^{1/2} \big) .
\end{equation}
Finally, an integration by parts, and \eqref{oswald} and \eqref{normal}, yield
\begin{equation}
\begin{aligned}
   (\div&(\bfs-\bft), \Ioh v)=-    ( \bfs-\bft , \nabla \Ioh v)\leq  \Vert \bfs-\bft\Vert_0 \Vert \nabla \Ioh v\Vert _0
   \\
   & \lesssim 
   \Vert \bfs-\bft\Vert_0 \Vert  v\Vert _ \etth
     \lesssim \mu^{-1/2}  \Vert \bfs-\bft\Vert_0.
   \end{aligned}
\end{equation} 
 Collecting  the estimates proves the claim.
\end{proof}

As noted before, above estimates give that the convergence rates are
 ${\mathcal O} (h^k) $ for CGG, AFW and \NEW, and  and ${\mathcal O}
 (h^{k+1}) $ for SGG, for a smooth solution.

 \section{Postprocessing of the displacement} \label{sec::postprocessing}

Since the work of Arnold and Brezzi \cite{AB} it is known that the solution can be postprocessed to yield an improved displacement approximation. Here we will use the postprocessing technique given in \cite{PPRAIRO,PPNM}. The postprocessing is done in two steps yielding an approximation in the spaces 
\begin{align} 
  V_h^* &=  \{\, v\in \Vtwo\, \vert \ v\vert _K\in [P_{l}(K)]^\dim\ \forall K\in \Ch\,\}, \label{eq::globaldisplpostspace}\\
  V_h^a &=V^*_h \cap \VV \label{eq::globaldisplaverspace},
\end{align}
with the choice $l=k$ for CGG, AFW and \NEW, and $l=k+1$ for SGG.

Further let $P_h^*: L^2(\Omega) \to V_h^*$  denote the $L^2$
projection on $V_h^*$.
 
{\em Postprocessing. Step I:~} The first step gives a  {\em
discontinuous} displacement by:  find $u_h^* \in V_h^*$ such that 
\begin{equation}\label{pp}
\begin{aligned}
P_h u_h^ * &= u_h,
\\  
(\nabla u_h^*,\nabla v)_K &= (\comp \bfs_h+\rho_h, \nabla v)_K \quad \forall v\in   (I-P_h) V_h^*\vert _K.
\end{aligned}
\end{equation}
 For this approximation we have the following estimate. We omit the proof and refer to the 
  analogues one in \cite{LS1}.
\begin{lemma}\label{discest}
It holds that
\begin{equation}
\norm{u-u_h^* }_ \etth \lesssim \Vert u-P_h^* u \Vert_ \etth+ \Vert \rho-\rho_h\Vert_0
+ \mu^{-1}\Big(\Vert\bfs-\bfs_h\Vert_{0} +\big(  \sum_{K\in \Ch} h_K^2 \Vert f- f_h\Vert_{0,K}^2\big)^{1/2} \Big).
\end{equation}
\end{lemma}

  {\em Postpostprocessing. Step II: ~} The second step gives a  continuous displacement approximation (used for the
  hypercircle technique below) by applying an averaging operator $\Ioh
  :V_h^* \to V_h^a$. Now let  $\uoh = \Ioh u_h^*$, then we have the
  following error estimate.
\begin{theorem} \label{Theo4}  It holds that
\begin{equation}
\begin{aligned}
\mu^{1/2} \norm{\varepsilon(u-&u_h^a) }_{0} 
  \lesssim \mu^{1/2} \Vert u-P_h^* u \Vert_ \etth+\mu^{1/2} \inf_{\eta \in R_h}\Vert \rho-\eta \Vert_0
\\
&+\mu^{-1/2} \big(\inf_{\tau \in \Sh^g }\Vert \bfs-\tau \Vert_0 +\big(  \sum_{K\in \Ch} h_K^2 \Vert f- f_h\Vert_{0,K}^2     \big)^{1/2}\big)
.
\end{aligned}
\end{equation}
\end{theorem}
For the proof we again refer to \cite{LS1}.

 \section{A posteriori error estimates by the hypercircle theorem}
 \label{sec::hypercircle} Of the previous  a posteriori estimates
 given in the literature,  the ones given in \cite{LV,CCDG, CCetal}
 all contain unknown constants. Only the work of Kim \cite{kim} gives
 an estimator with known constants. Below we will use the hypercircle theorem to derive  an alternative estimator with  known constants. The proof of this classical result is given in
 \cite{NecasH} and recalled in \cite{LS1}.

\begin{theorem}\label{hypc} {\rm (The Prager-Synge hypercircle theorem)}   Suppose that:

\begin{itemize}
\item The stress $\Sigma\in \Hdiv\cap \SY $ is symmetric $\Sigma=\Sigma^T, $
      statically admissible $\div \Sigma +f=0 \mbox{ in } \Omega $, and $\Sigma n=g \mbox{ on } \Gamma_N$.
\item The displacement $U\in [H^1(\Omega)]^\dim$ is kinematically admissible; $U\vert_{\Gamma_D} =0. $
\end{itemize}
  Then it holds
    \begin{equation}
\begin{aligned}
    \enorm{\bfs -\Sigma} ^2+ \enorm{\bfs -\elas \varepsilon(U)}^2   =   \enorm{\Sigma-\elas \varepsilon(U)}^2,
    \end{aligned}
    \end{equation}
   and
     \begin{equation}
    \begin{aligned}
    \enorm{\bfs -&\frac{1}{2} \big(\Sigma+ \elas \varepsilon(U)\big)}       = \frac{1}{2}  \enorm{\Sigma-\elas \varepsilon(U)}.
  \end{aligned}
      \end{equation}
      \end{theorem}
     
  Of two reasons  this theorem cannot be applied directly. The finite
  element solution is not exactly symmetric and, in general, the
  equilibrium equations are not exactly valid. Therefore we consider
  two auxiliary problems. The first one is the continuous problem with the
  loading in the finite element spaces. To this end, let
  $P_K=P_h\vert_K$, where $P_h$ is the projection \eqref{projpro}. In
  \cite{LS1}  we proved the following estimate.
  \begin{lemma}\label{pert1}
  For the solution  $(\bar \bfs,\bar \rho, \bar u) \in \E$  to 
\begin{equation}
\B (\bar \bfs, \bar \rho, \bar u; \bft, \eta, v) =(P_hf,v) + \langle  Q_h g, v \rangle_{\Gamma_N}  \quad \forall(\bft,\eta, v) \in \E,
\end{equation}
it holds that 
\begin{equation}
\norm{(\bfs- \bar \bfs, \rho-\bar \rho, u-\bar u) }_\E \leq C_1 osc(f)+ C_2 osc(g),
\end{equation}
  with  constants $C_1, C_2 > 0$,
  respectively, and
  \begin{equation}  osc(f)=\big(  \sum_{K\in \Ch} h_K^2 \Vert f-P_K f \Vert_{0,K}^2     \big)^{1/2},
  \end{equation}
  and 
  \begin{equation}
  osc(g)=\big(  \sum_{E\subset \Gamma_N} h_E \Vert  g-Q_E g   \Vert_{0,E}^2\big)^{1/2}.
  \end{equation}
  \end{lemma}
  Next, we construct a symmetric approximation of $\bfs_h$. 
 
 
   Here $\korn$  denotes the constant in the Korn inequality
  \begin{equation}\label{korn}
  \Vert \omega(v) \Vert_0 \leq \korn   \Vert \varepsilon(v) \Vert_0  \quad \forall v\in \VV.
  \end{equation}
 \begin{lemma}\label{pert2}
For the solution  $(\st,\rt, \zt) \in \E$ of
 \begin{align}\label{auxprob}
 \B(\st,\rt, \zt;\bft, \eta , v)=
 ( \mathcal{C} \bfs_h, \bft) 
 +(P_h f,v) 
 + \langle Q_h  g, v \rangle_{\Gamma_N}  
 & 
 \quad \forall (\bft,\eta,v)\in \E , 
 \end{align}
    it holds that
    \begin{equation}\label{auxequi}
    \div \st=-P_h f \ \mbox{ in } \Omega, \quad \st n=Q_h g  \ \mbox{ on } \Gamma_N,
    \end{equation}
    and
\begin{equation} \enorm{ \bfs_h -\st } \leq \frac{{\sqrt{5}+1}}{4}  \big(C_\Omega+1\big) \enorm{\bfs_h-\bfs_h^T},
\end{equation}
\begin{proof} 
Testing \eqref{auxprob} with $v$ yields
\begin{equation}
(\st, \nabla v)=(P_h f,v) + \langle Q_h  g, v \rangle_{\Gamma_N} ,
\end{equation}
which gives \eqref{auxequi}.

It holds
\begin{equation}
\begin{aligned}
  \B(\st&-\bfs_h ,\rt, \zt;\bft, \eta, v) 
 = \B(\st ,\rt, \zt;\bft, \eta, v)- \B( \bfs_h ,0, 0;\bft, \eta, v)
 \\&=  ( \mathcal{C} \bfs_h, \bft) +(P_h f,v) + \langle Q_h  g, v \rangle_{\Gamma_N} 
 \\&\qquad  - ( \mathcal{C} \bfs_h, \bft) + (\bfs_h,\eta) - (\bfs_h, \nabla v) 
 \\& =   (P_h f,v) + \langle Q_h  g, v \rangle_{\Gamma_N}   + (\bfs_h,\eta) - (\bfs_h, \nabla v) .
 \end{aligned}
 \end{equation}
Since
\begin{equation} 
  -(\bfs_h, \nabla v) = (\div \bfs_h, v) - \langle \bfs_h n,v \rangle_{\Gamma_N} = -(P_h f,v) - \langle Q_h  g, v \rangle_{\Gamma_N} ,
  \end{equation}
 we obtain
 \begin{equation}
   \B(\st-\bfs_h ,\rt, \zt;\bft, \eta, v) =(\bfs_h,\eta) .
 \end{equation} 
 Since $\eta$ is skew symmetric it holds
 \begin{equation}
 (\bfs_h,\eta)=\frac{1}{2} (\bfs_h-\bfs_h^T,\eta).
 \end{equation}
By adding and subtracting $\omega(v)$
 \begin{equation} 
 \B(\st- \bfs_h ,\rt, \zt;\bft, \eta, v) 
  = \frac{1}{2} (\bfs_h-\bfs_h^T, \eta - \omega(v)) +  \frac{1}{2} (\bfs_h -\bfs_h^T, \omega(v)).
 \end{equation}
 Theorem \ref{energystab} and the inequality \eqref{korn} then yields

 \begin{equation}
 \Big(\frac{\sqrt{5}-1}{2} \Big) \enorm{\st- \bfs_h } \leq \frac{1}{2}\big( 1+C_\Omega \big)  \enorm{\bfs_h -\bfs_h^T },
 \end{equation}
 which proves the claim.
\end{proof}
 \end{lemma}

The hypercircle  estimate is now.
\begin{theorem} \label{th::estimator}
It holds that
 \begin{equation} \label{avarageest}
    \begin{aligned}
     \enorm{\bfs -\frac{1}{2} \big(\bfs_h&+ \elas \varepsilon(\uoh)\big)}       
 \leq\frac{1}{2}  \enorm{\bfs_h-\elas \varepsilon(\uoh)}  
 \\ &
  +  \frac{{\sqrt{5}+1}}{4}  \big(C_\Omega+1\big) \enorm{\bfs_h-\bfs_h^T}    + C_1 osc(f)+ C_2 osc(g),
  \end{aligned}
  \end{equation}
  and
  \begin{equation}  \label{sumest}
    \begin{aligned}
    \enorm{\bfs -\bfs_h}   +    \enorm{\bfs - \elas \varepsilon(\uoh)}     
      \lesssim   \enorm{\bfs_h-\elas \varepsilon(\uoh)}  
    &  +   \enorm{\bfs_h-\bfs_h^T}  \\ &  +   osc(f)+   osc(g). 
  \end{aligned}
    \end{equation}
     
    \end{theorem}
\begin{proof}   
Applying the Hypercircle theorem with the symmetric stresses $\bar \bfs$,  $\tilde \bfs$, and $\uoh$, gives 
\begin{equation}
  \enorm{\bar \bfs -\frac{1}{2} \big(\tilde \bfs+ \elas \varepsilon(\uoh)\big)}       
    \leq\frac{1}{2}  \enorm{\tilde \bfs -\elas \varepsilon(\uoh)}  .
\end{equation}
Using the triangle inequality we get
\begin{equation}
\enorm{  \bfs -\frac{1}{2} \big(  \bfs_h+ \elas \varepsilon(\uoh)\big)}  \leq 
 \enorm{\bar \bfs -\frac{1}{2} \big(\tilde \bfs+ \elas \varepsilon(\uoh)\big)} + \enorm{\bfs -\bar \bfs } +\frac{1}{2} \enorm{\tilde\bfs -\bfs_h },
\end{equation}
and
\begin{equation}
\enorm{\tilde \bfs -\elas \varepsilon(\uoh)} \leq \enorm{\bfs_h -\elas \varepsilon(\uoh)}  + \enorm{\tilde \bfs -\bfs_h} .
\end{equation}
Combining gives
\begin{equation}
\enorm{  \bfs -\frac{1}{2} \big(  \bfs_h+ \elas \varepsilon(\uoh)\big)}  \leq 
\enorm{\bfs_h -\elas \varepsilon(\uoh)}  +\enorm{\bfs -\bar \bfs } + \enorm{\tilde\bfs -\bfs_h },
\end{equation}
and \eqref{avarageest} follows  from Lemma \ref{pert1} and Lemma \ref{pert2}. 

In the same manner \eqref{sumest} follows from the first Hypercircle identity.
\end{proof}
\begin{remark} For $f$   smooth, it holds that
\begin{equation}
osc(f) = \mathcal{O}(h^{k+1}.
\end{equation}
which is a higher order term for all methods but the CGG. 
However, in most real
engineering problems the loading $f$ is a constant and $osc(f)$ vanishes.

For a smooth traction $g$ we have  
\begin{equation}
 osc(g) = \begin{cases}\mathcal{O}(h^{k+3/2}) \quad \mbox{for SGG and AFW},
\\ \mathcal{O}(h^{k+1/2}) \quad \mbox{for CGG and \NEW}.\end{cases}
\end{equation}
Hence,   $osc(g)$ is a higher order term for all methods. 

For the modification \NEW~ of the AFW method there is the option not to
 reduce the polynomial order for the edges on $\Gamma_N$ which would
 lead to 
 $ \mathcal{O}(h^{k+3/2}) $, i.e. an convergence of order $\mathcal{O}(h^{3/2}) $  faster than that of the convengence order $\mathcal{O}(h^{k})$.
\end{remark}

  \section{An a   posteriori estimator uniformly valid in the incompressible limit}
  
  The drawback of the estimate by  the hypercircle argument is that it
 is formulated in terms of $\enorm{\cdot}$ which, unfortunately,
 ceases to be a norm in the incompressible limit $\lambda\to \infty$
 and that the stress computed from the displacement, i.e.
   \begin{equation}  \label{compstress} 
 \elas \varepsilon(u_h^a)  =   2\mu \varepsilon(u_h^a)  +\lambda \div u_h^a I ,
 \end{equation} 
 grows without limit unless $\div u_h^a $ will vanish identically in the limit. For two space dimensions it
 is well known \cite{MR1174468,MR813691} that this requires piecewise
 polynomials of degree four or higher.


In this section we will therefore derive the following  a posteriori estimate uniformly valid with respect to the second Lam\'e parameter.
  \begin{theorem} \label{th::estimatorinc}
   It holds that
  \begin{equation}
  \begin{aligned}
   \Norm{(\bfs-\bfs_h , \rho-\rho_h, u-\uoh)}
 & \lesssim
 \mu^{1/2}  \Vert \comp \bfs_h+\rho_h -\nabla \uoh   \Vert_0 +\mu^{-1/2} \norm{\bfs_h -\bfs_h^T}_0
  \\
   &\qquad +\mu^{-1/2}  \Big( osc(f)+osc(g) \Big).
  \end{aligned}
  \end{equation}
  \end{theorem}
\begin{proof}
 By Theorem \ref{incompstab} there exists  $(\bft,\eta, v) \in \E $, with    $\Norm{(\bft,\eta, v)} =1$, such that
  \begin{equation}
  \begin{aligned}
 &   \Norm{(\bfs-\bfs_h , \rho-\rho_h, u-\uoh)} \lesssim
  \B( \bfs-\bfs_h , \rho-\rho_h,  u-\uoh ; \bft, \eta, v)
  \\&
  =(\comp (\bfs-\bfs_h),\bft ) +(\rho-\rho_h, \tau) - (\nabla(u-\uoh),\bft) 
  \\
  & \quad +\big((\bfs-\bfs_h)-(\bfs^T -\bfs_h^T), \eta\big) - (\nabla v, \bfs-\bfs_h).
  \end{aligned} 
\end{equation}
Since $ \comp \bfs +\rho -\nabla u=0 $
 we have
  \begin{equation}
  \begin{aligned}
    (\comp (\bfs-\bfs_h),\bft ) &+(\rho-\rho_h, \tau) - (\nabla(u-\uoh),\bft) =  - (\comp \bfs_h +\rho_h -\nabla u_h^a,\bft )  
    \\
 &\leq \Vert  \comp \bfs_h +\rho_h -\nabla u_h^a\Vert_0\Vert \bft\Vert_0
 \\& \leq \mu^{1/2} \Vert  \comp \bfs_h +\rho_h -\nabla u_h^a\Vert_0.
 \end{aligned}
 \end{equation}
 From  $\bfs-\bfs^T=0$ it follows that 
 \begin{equation}
\big((\bfs-\bfs_h)-(\bfs^T -\bfs_h^T), \eta\big) = (\bfs_h^T-\bfs_h, \eta) \leq \norm{\bfs_h^T-\bfs_h}_0  \norm{  \eta}_0 \lesssim \mu^{-1/2}  \norm{\bfs_h^T-\bfs_h}_0 .
 \end{equation}
   Since $\div \bfs_h = - P_h f$ and $\sigma_h n=Q_hg$ we get 
  \begin{equation}
  \begin{aligned}
  - (&\nabla v, \bfs-\bfs_h) = (\div(\bfs-\bfs_h), v  ) - \langle  (\bfs-\bfs_h)n, v \rangle_{\Gamma_N}
  \\& = (P_h f -f , v  ) + \langle  Q_h g - g, v \rangle_{\Gamma_N}  .
  \end{aligned}
  \end{equation}
  By using the properties of the projection and Korn's inequality one gets (see \cite[Theorem 6]{LS1})
   \begin{equation}
  \begin{aligned}
    (P_h f -f  , v  ) 
             & \leq C\big( \sum_{K\in \Ch} h_K^2 \norm{f-P_h f }_{0,K}^2\big)^{1/2}  \norm{\varepsilon( v )  }_{0}
                         \\&     \leq \mu^{-1/2}C\big( \sum_{K\in \Ch} h_K^2 \norm{f-P_h f }_{0,K}^2\big)^{1/2}.
  \end{aligned}
  \end{equation}
  Using the trace theorem we similarly obtain
  \begin{equation}
\langle  Q_h g - g, v\rangle_{\Gamma_N} \leq \mu^{-1/2} C \big( \sum_{E\subset \Gamma_N} h_E \Vert g-Q_K g \Vert_{0,E} ^2 \big)^{1/2} . 
  \end{equation}
  Hence, we have 
  \begin{equation}
   - (\nabla v, \bfs-\bfs_h) \lesssim \mu^{-1/2}  \Big( osc(f)+osc(g) \Big).
   \end{equation}
   The assertion follows by collecting the above estimates.
   \end{proof}


   \section{Numerical examples}
   In this section we consider several numerical examples to strengthen
   our theoretical findings. All examples were implemented in the finite
   element library Netgen/NGSolve, see \cite{netgen}. In our
   computations, we take into account the relative errors in strain
   energy for the stress directly obtained from the method, and computed
   from the post processed displacement $u_h^a$, and the relative error
   of Theorem~\ref{th::estimator},
   \begin{equation*} 
     e_\mathcal{C}^\sigma  = \frac{\|
   \sigma - \sigma_h \|_\mathcal{C} }{\| \sigma\|_\mathcal{C}},
   \quad \quad 
   e_\mathcal{C}^{u}  = \frac{\| \sigma -
   \mathcal{A}\eps(u_h^a) \|_\mathcal{C} }{\| \sigma\|_\mathcal{C}},
   \quad \quad
   {e}^{\operatorname{mean}}_\mathcal{C}   = \frac{
       \| \sigma - \frac{1}{2} (\sigma_h + \mathcal{A}\eps(u_h^a)) \|_\mathcal{C} 
     }{\|
     \sigma\|_\mathcal{C}}.
   \end{equation*}
   In addition to that we make use of the (relative) estimators
   \begin{align*}
     \eta &= \frac{\frac{1}{2} \| \sigma_h -  \mathcal{A}\eps(u_h^a) \|_\mathcal{C} 
     +  \|\sigma_h - \sigma_h^T \|_\mathcal{C}
     }{\|
     \sigma\|_\mathcal{C}}, \\
     \eta^{\operatorname{inc}} &= \frac{
     \mu^{1/2} \| \mathcal{C}\sigma_h + \rho_h -
     \nabla(u_h^a)\|_0 + \mu^{-1/2} \|\sigma_h - \sigma_h^T \|_0
     }{\mu^{-1/2}\| \sigma\|_0},
   \end{align*}

   and the scaled errors from Theorem~\ref{th::estimatorinc} 
   
   

\begin{align*}
     e_{\operatorname{inc}}^{\sigma} =\frac{  \|
   \sigma - \sigma_h\|_0 }{ \| \sigma\|_0},
   \quad \quad 
     e_{\operatorname{inc}}^{u} =\frac{ \|
     \eps(u) - \eps(u_h^a)\|_0 }{ \|   \eps(u) \|_0}.
   \end{align*}
   
   If the exact solution $\sigma$ is not known, the scaling with $\|
   \sigma \|_\mathcal{C}$ and $\| \sigma \|_0$ is replaced by a
   scaling with $\|\sigma_h \|_\mathcal{C}$ and $\| \sigma_h \|_0$
   where $\sigma_h$ is the discrete solution from the finest mesh.
   
   In Table~\ref{tab::disc} we have listed the
   finite element families that we consider for our computations and
   give the convergence order of the above listed relative errors
   for smooth solutions in terms of the number of elements $\Nel$. Note that
   in two dimensions the families from \cite{family} and \cite{GG}
   coincide and only differ in three dimensions. To simplify the notation
   we always use the abbreviation SGG but used the family from \cite{GG}
   for the three dimensional example.
   
   For all examples we use an adaptive refinement with the estimator
   from Theorem~\ref{th::estimator} and for the incompressible limit
   the estimator given in Theorem~\ref{th::estimatorinc} (see 
   $\eta, \eta^{\operatorname{inc}}$ above). The adaptive mesh
   refinement loop is defined as usual by 
   \begin{align*}
      \mathrm{SOLVE} \rightarrow \mathrm{ESTIMATE} \rightarrow \mathrm{MARK} \rightarrow \mathrm{REFINE} \rightarrow \mathrm{SOLVE} \rightarrow \ldots
   \end{align*}
   During the marking routine we mark an element  $K \in \Ch$ for
   refinement if $\eta(K) \geq \frac{1}{4} \max\limits_{K \in \Ch}
   \eta(K)$ or $\eta^{\operatorname{inc}}(K) \geq \frac{1}{4}
   \max\limits_{K \in \Ch} \eta^{\operatorname{inc}}(K)$, where $\eta(K)$
   and $\eta^{\operatorname{inc}}(K)$ are the local contributions of the
   corresponding estimators of Theorem~\ref{th::estimator} and
   Theorem~\ref{th::estimatorinc}, respectively. After that, the mesh
   refinement algorithm refines the marked elements plus further elements
   to guarantee a regular triangulation.

   \begin{table}
     \begin{center}
       \begin{tabular}{lcc}
         method & abbr. & conv.  \\
         \midrule
         Stenberg, Gopalakrishnan, Guzm\'an family \cite{family, GG}, $k=2$ & SGG2 & $\Nel^{-3/d}$ \\
         Stenberg, Gopalakrishnan, Guzm\'an family \cite{family, GG}, $k=3$ & SGG3 & $\Nel^{-4/d}$ \\
         Arnold-Falk-Winther family \cite{AFW}, $k=2$ & AFW & $\Nel^{-2/d}$ \\
         The method from \eqref{sspace_rafw}, $k=2$ & \NEW & $\Nel^{-2/d}$ \\
       \end{tabular}
       \vspace{3mm}
       \caption{Considered discretizations for our computations.} \label{tab::disc}
     \end{center}
   \end{table}
   
   Below we always choose the Young modulus $E$ and the Poisson ratio
   $\nu$ rather than the Lam\'e parameters. They are related to each
   other by
   \begin{align*}
     \lambda = \frac{E \nu}{(1 + \nu)(1 - 2\nu)}, \quad \textrm{and} \quad \mu = \frac{E}{2(1+\nu)}.
   \end{align*}
   
   For the two dimensional problem we always consider the plane strain
   setting.  
   \subsection{L-shape example} \label{ex::lshape}
   
   We solve the L-shape example from Section 10.3.2 of
   \cite{szabo1991finite} with an adaptive mesh refinement. The domain
   $\Omega$ is 
   \begin{align*}
   \Omega = \{ (x,y): |x| + |y| \le {2}^{1/2} l \} \setminus  \{ (x,y): |x- 2^{-1/2} l| + |y| \le 2^{-1/2} l \}.
   \end{align*}
   We choose $l = 1$, $E = 1$ and $\nu=0.3$ or $\nu = 0.4999$.  The exact
   displacement field, up to rigid-body displacements and rotations, is
   given by 
   \begin{align*}
     u_x &= \frac{1}{2\mu} r^\alpha  ( (\kappa - Q(\alpha +1))  \cos(\alpha \theta) - \alpha \cos((\alpha - 2)\theta) ), \\
     u_y &= \frac{1}{2\mu} r^\alpha  ( (\kappa + Q(\alpha +1))  \sin(\alpha \theta) + \alpha \sin((\alpha - 2)\theta) ),
   \end{align*}
   and the stress components are
   \begin{align*}
     \sigma_{xx} &= \alpha  r^{\alpha-1}  ( (2 - Q(\alpha + 1))  \cos((\alpha-1)  \theta) - (\alpha-1) \cos((\alpha - 3)\theta) ), \\
     \sigma_{yy} &= \alpha  r^{\alpha-1}  ( (2 + Q(\alpha + 1))  \cos((\alpha-1)  \theta) + (\alpha-1) \cos((\alpha - 3)\theta) ), \\ 
     \sigma_{xy} &= \alpha  r^{\alpha-1}  ( (\alpha-1)  \sin((\alpha -3)\theta) + Q (\alpha+1)  \sin((\alpha-1)\theta)). 
   \end{align*}
   where  $\alpha = 0.544483737 $ and $Q = 0.543075579$ are fixed
   constants.
   
\begin{figure}
     \begin{center}
       \resizebox{.95\linewidth}{!}{\input{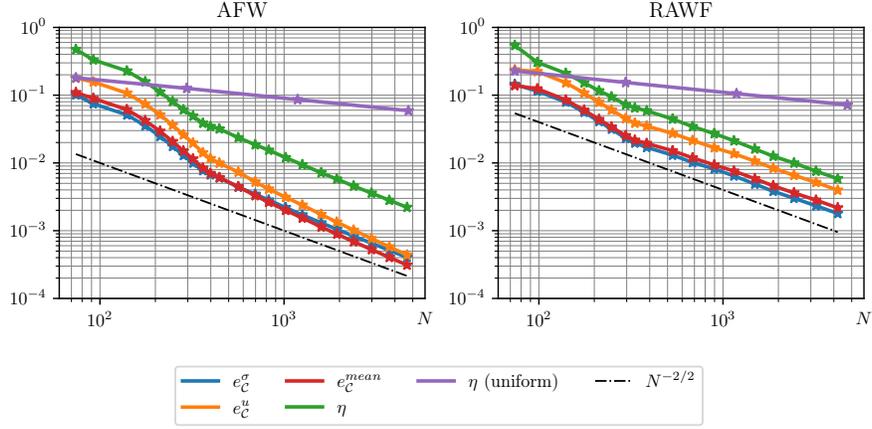}}
     \end{center}
     \caption{Convergence of the AFW and \NEW~ element for the L-shape example for $\nu = 0.3$.}
     \label{fig::AFWlshape}
   \end{figure}

   In Figure~\ref{fig::AFWlshape} and Figure~\ref{fig::STENlshape} we
   present the convergence of the errors $e_{\mathcal{C}}^\sigma$,
   $e_{\mathcal{C}}^u$ and $e_{\mathcal{C}}^{\operatorname{mean}}$ and the estimator $\eta$ for a moderate Poisson
   ratio $\nu = 0.3$. 
   From the figures we see
   that all errors converge with optimal order. To demonstrate the
   drastic decrease of the error for an adaptive mesh refinement
   strategy, we included the estimator $\eta$ also for a uniform
   refinement. Since the exact solution is in the Sobolev space $H^s$,
   with $s < 1.54$, a uniform mesh only yields a convergence rate of
   $\mathcal{O}(h^{0.54})$, i.e. $\mathcal{O}(\Nel^{-0.27})$. 
   
   In Figure~\ref{fig::AFWlshape_inc} and
   Figure~\ref{fig::STENlshape_inc} we solve the same problem but with
   a Poisson ratio $\nu = 0.4999$ and using the estimator for the
   incompressible setting from Theorem~\ref{th::estimatorinc}. Again
   we plot the errors $e_{\mathcal{C}}^\sigma$ and $e_{\mathcal{C}}^u$
   but also include $e_{\operatorname{inc}}^\sigma$,
   $e_{\operatorname{inc}}^u$ and the relative estimator
   $\eta^{\operatorname{inc}}$. 
   The errors  $e_{\mathcal{C}}^u$,
   $e_{\mathcal{C}}^{\operatorname{mean}}$ and $\eta$,  deteriorate
   for values of $\nu$ near $1/2$, and are therefore not plotted.
   Thus, in practice
   $\eta$ should not be used for this setting. On the other hand,
   $e_{\mathcal{C}}^\sigma, e_{\operatorname{inc}}^\sigma$ and $e_{\operatorname{inc}}^u$ as well
   as  $\eta^{\operatorname{inc}}$ are not affected by the choice of
   $\nu$ and are much smaller.

   \begin{figure}
     \begin{center}
       \resizebox{.95\linewidth}{!}{\input{graphics/lshape_STEN23.pgf}}
     \end{center}
     \caption{Convergence of the SGG2 and SGG3 element for the L-shape example for $\nu = 0.3$.}
     \label{fig::STENlshape}
   \end{figure}
   
   \begin{figure}
     \begin{center}
       \resizebox{.95\linewidth}{!}{\input{graphics/lshape_AFWRAFW_inc.pgf}}
     \end{center}
     \caption{Convergence of the AFW and \NEW~ element for the L-shape example for $\nu = 0.4999$.}
     \label{fig::AFWlshape_inc}
   \end{figure}

   \begin{figure}
     \begin{center}
       \resizebox{.95\linewidth}{!}{\input{graphics/lshape_STEN23_inc.pgf}}
     \end{center}
     \caption{Convergence of the SGG2 and SGG3 element for the L-shape example for $\nu = 0.4999$.}
     \label{fig::STENlshape_inc}
   \end{figure}

   \subsection{Cook's membrane} \label{ex::cook} 
\begin{figure}
     \begin{center}
       \begin{minipage}{0.4\textwidth}
         \includegraphics[width=\textwidth]{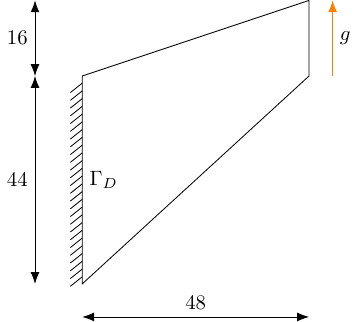}
       \end{minipage}
       \begin{minipage}{0.4\textwidth}
         \begin{tikzpicture}
           \node[anchor=south west] (bar1) at (4.8,0) {
               \pgfplotscolorbardrawstandalone[colormap/jet
               ,colorbar sampled, point meta min=0, point meta max=0.5, 
               colorbar style={samples=16, width=0.2cm, height = 4.5cm, ytick={0, 0.1, 0.2, 0.3, 0.4, 0.5}, 
               scaled y ticks=false, 
               yticklabel style={style={font=\footnotesize}, /pgf/number format/fixed,  /pgf/number format/precision=3}}]};
           \node[anchor=south west] (numex1) at (0,0)  {\includegraphics[width=\textwidth]{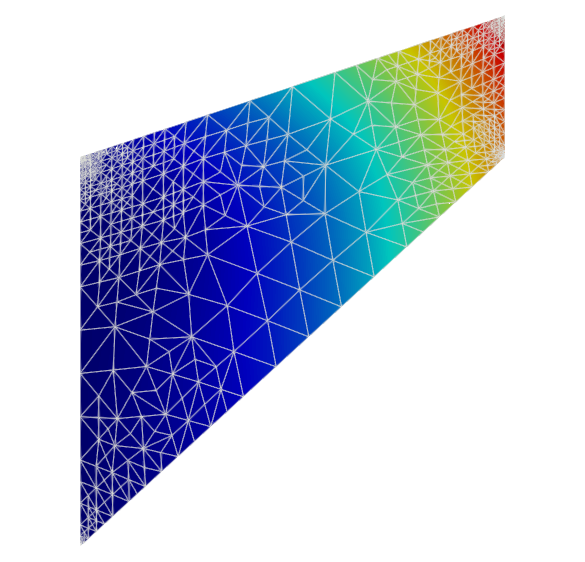}};
         \end{tikzpicture}
         
       \end{minipage}
     \end{center}
       \caption{Computational domain and absolute value of the
       displacement on the adaptively refined triangulation of
       example~\ref{ex::cook} using the SGG2 element.}\label{fig::cook}
     \end{figure}

We consider the Cook's
   membrane benchmark problem. The geometry $\Omega$ is depicted in the
   left picture of Figure~\ref{fig::cook}. The problem describes a
   tapered beam which is clamped on th left boundary $\Gamma_D$. On the
   opposite boundary we apply a traction force $g = (0,10^{-3})$ and we
   choose $E = 1$ and a Poisson ratio $\nu = 0.3$. All other boundary
   edges are free (i.e. $g = 0$). The solution is known to have singularities
   in the corners, i.e. an adaptive refinement is desirable. In the right
   picture of Figure~\ref{fig::cook} we have drawn the absolute value
   of the displacement and the adaptively refined triangulation. One can
   clearly see that the estimator successfully detects the corner
   singularities. In Figure~\ref{fig::cookconv} we have plotted the error
   estimator $\eta$ for the elements from Table~\ref{tab::disc}. We
   observe that all elements provide an optimal convergence. Again, to
   show the drastic decrease of the error we have plottet the estimator
   $\eta$ for the \NEW ~element also for a uniform refinement.

   \begin{figure}
     \begin{center}
       \resizebox{.7\linewidth}{!}{\input{graphics/cook2.pgf}}
     \end{center}
     \caption{Convergence of the error estimator $\eta$ for example~\ref{ex::cook}.}\label{fig::cookconv}
   \end{figure}

     \subsection{Clamped Fichera corner} \label{ex::fichera}
   
\begin{figure}
       \begin{center}
           \begin{tikzpicture}
             \node[anchor=south west] (bar1) at (9,1) {
                 \pgfplotscolorbardrawstandalone[colormap/jet
                 ,colorbar sampled,point meta min=0, point meta max=0.02, 
                 colorbar style={samples=17, width=0.2cm, height = 6cm, ytick={0, 0.005, 0.01, 0.015, 0.02}, 
                 scaled y ticks=false, 
                 yticklabel style={style={font=\footnotesize}, /pgf/number format/fixed,  /pgf/number format/precision=3}}]};
             \node[anchor=south west] (numex1) at (0,0)  {\includegraphics[width=0.7\textwidth]{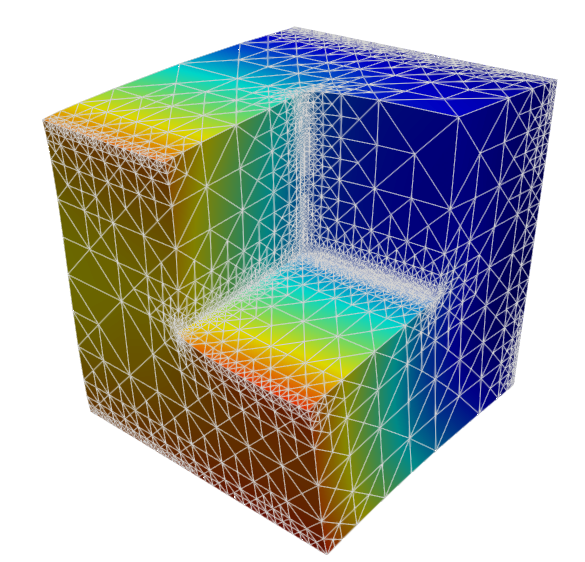}};
           \end{tikzpicture}
       \end{center}
         \caption{Absolute value of the displacement on the deformed
         adaptively refined triangulation of example \ref{ex::fichera}
         using the SGG2 element.}\label{fig::fichera}
       \end{figure}

\begin{figure}
         \begin{center}
           \resizebox{.7\linewidth}{!}{\input{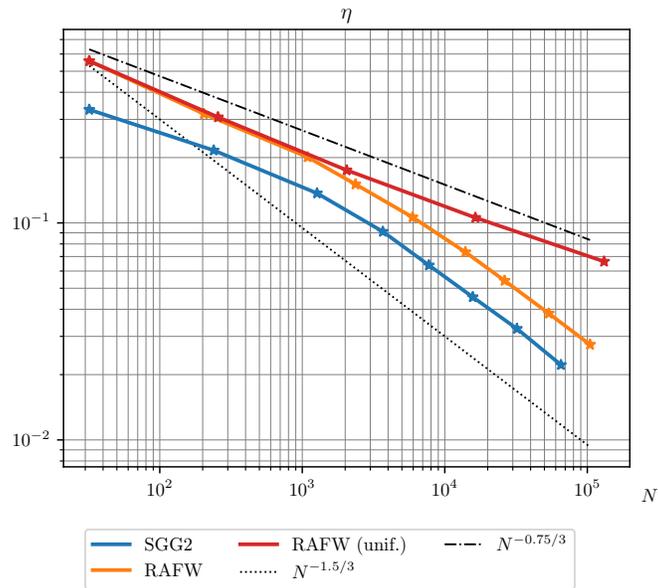}}
         \end{center}
         \caption{Convergence of the error estimator $\eta$ for example \ref{ex::fichera}}\label{fig::ficheraconv}
       \end{figure}

     Wo consider a clamped Fichera corner example. The geometry is given
     by $\Omega = (-1,1)^3 \setminus (0,1)^3$ and we choose the
     parameters $E = 1$ and $\nu = 0.3$. We consider the object to be
     clamped at $x = -1$ and apply a traction force $g = (0,0,-10^{-3})$
     on the opposite side where $x=1$. All other boundaries are free
     (i.e. $g = 0$). In contrast to the two dimensional setting the
     solution now includes vertex and edge singularities, i.e. even for
     an adaptively refined triangulation we might not see the "optimal"
     convergence if we do not apply an anisotropic refinement strategy,
     see \cite{MR1716824}. Since this is out of scope of this work we
     applied the same adaptive refinement as for the other examples. In
     Figure~\ref{fig::fichera} we have plotted the absolute value of the
     displacement using the SGG2 finite element and deformed the geometry
     with a scaling factor of 20 since otherwise the deformation would
     not be visible. In addition we show the adaptively refined
     triangulation (from the last step). One can clearly see that the
     estimator successfully detects the edge and vertex singularities
     i.e. that the corresponding elements are marked and refined during
     the computation loop. In Figure~\ref{fig::ficheraconv} we plot the
     estimator $\eta$ for the SGG2 and the \NEW ~element using an adaptive
     refinement and included also the convergence using a uniform
     refinement for the \NEW ~element. Although we do not see the optimal
     convergence (as noted above) we still see the drastic decrease in
     the error when using an adaptive mesh refinement.

\bibliography{stenberg-2}

\begin{thebibliography}{10}

\bibitem{AT}
M.~Amara and J.~M. Thomas.
\newblock Equilibrium finite elements for the linear elastic problem.
\newblock {\em Numer. Math.}, 33(4):367--383, 1979.

\bibitem{MR1716824}
Thomas Apel.
\newblock {\em Anisotropic finite elements: local estimates and applications}.
\newblock Advances in Numerical Mathematics. B. G. Teubner, Stuttgart, 1999.

\bibitem{AB}
D.~N. Arnold and F.~Brezzi.
\newblock Mixed and nonconforming finite element methods: implementation,
  postprocessing and error estimates.
\newblock {\em RAIRO Mod\'el. Math. Anal. Num\'er.}, 19(1):7--32, 1985.

\bibitem{AAW3D}
Douglas~N. Arnold, Gerard Awanou, and Ragnar Winther.
\newblock Finite elements for symmetric tensors in three dimensions.
\newblock {\em Math. Comp.}, 77(263):1229--1251, 2008.

\bibitem{PEERS}
Douglas~N. Arnold, Franco Brezzi, and Jim Douglas, Jr.
\newblock P{EERS}: a new mixed finite element for plane elasticity.
\newblock {\em Japan J. Appl. Math.}, 1(2):347--367, 1984.

\bibitem{ADG}
Douglas~N. Arnold, Jim Douglas, Jr., and Chaitan~P. Gupta.
\newblock A family of higher order mixed finite element methods for plane
  elasticity.
\newblock {\em Numer. Math.}, 45(1):1--22, 1984.

\bibitem{MR2249345}
Douglas~N. Arnold, Richard~S. Falk, and Ragnar Winther.
\newblock Differential complexes and stability of finite element methods. {II}.
  {T}he elasticity complex.
\newblock In {\em Compatible spatial discretizations}, volume 142 of {\em IMA
  Vol. Math. Appl.}, pages 47--67. Springer, New York, 2006.

\bibitem{AFW}
Douglas~N. Arnold, Richard~S. Falk, and Ragnar Winther.
\newblock Finite element exterior calculus, homological techniques, and
  applications.
\newblock {\em Acta Numer.}, 15:1--155, 2006.

\bibitem{AFW07}
Douglas~N. Arnold, Richard~S. Falk, and Ragnar Winther.
\newblock Mixed finite element methods for linear elasticity with weakly
  imposed symmetry.
\newblock {\em Math. Comp.}, 76(260):1699--1723 (electronic), 2007.

\bibitem{MR1174468}
Ivo Babu\v{s}ka and Manil Suri.
\newblock Locking effects in the finite element approximation of elasticity
  problems.
\newblock {\em Numer. Math.}, 62(4):439--463, 1992.

\bibitem{MR3097958}
Daniele Boffi, Franco Brezzi, and Michel Fortin.
\newblock {\em Mixed finite element methods and applications}, volume~44 of
  {\em Springer Series in Computational Mathematics}.
\newblock Springer, Heidelberg, 2013.

\bibitem{brennerkorn}
Susanne~C. Brenner.
\newblock Korn's inequalities for piecewise {$H^1$} vector fields.
\newblock {\em Math. Comp.}, 73(247):1067--1087, 2004.

\bibitem{CCetal}
C.~Carstensen, G.~Dolzmann, S.~A. Funken, and D.~S. Helm.
\newblock Locking-free adaptive mixed finite element methods in linear
  elasticity.
\newblock {\em Comput. Methods Appl. Mech. Engrg.}, 190(13-14):1701--1718,
  2000.

\bibitem{CCDG}
Carsten Carstensen and Georg Dolzmann.
\newblock A posteriori error estimates for mixed {FEM} in elasticity.
\newblock {\em Numer. Math.}, 81(2):187--209, 1998.

\bibitem{CGG}
Bernardo Cockburn, Jayadeep Gopalakrishnan, and Johnny Guzm{\'a}n.
\newblock A new elasticity element made for enforcing weak stress symmetry.
\newblock {\em Math. Comp.}, 79(271):1331--1349, 2010.

\bibitem{DiPietroErn}
Daniele~Antonio Di~Pietro and Alexandre Ern.
\newblock {\em Mathematical aspects of discontinuous {G}alerkin methods},
  volume~69 of {\em Math\'ematiques \& Applications (Berlin) [Mathematics \&
  Applications]}.
\newblock Springer, Heidelberg, 2012.

\bibitem{FdV2}
B.~M. Fraeijs~de Veubeke.
\newblock Discretization of rotational equilibrium in the finite element
  method.
\newblock In {\em Mathematical aspects of finite element methods ({P}roc.
  {C}onf., {C}onsiglio {N}az. delle {R}icerche ({C}.{N}.{R}.), {R}ome, 1975)},
  pages 87--112. Lecture Notes in Math., Vol. 606. Springer, Berlin, 1977.

\bibitem{FdV3}
B.~M. Fraeijs~de Veubeke and A.~Millard.
\newblock Discretization of stress fields in the finite element method.
\newblock {\em J. Franklin Inst.}, 302(5-6):389--412, 1976.
\newblock Basis of the finite element method.

\bibitem{FdV}
B.M. Fraejis~de Veubeke.
\newblock Displacement and equilibrium models in the finite element method.
\newblock In O.C Zienkiewics and G.S. Holister, editors, {\em Stress analysis},
  pages 145--197. Wiley, 1965.

\bibitem{FdV-W}
B.M. Fraejis~de Veubeke.
\newblock Stress function approach.
\newblock {\em Proc. of the World Congress on Finite Element Methods in
  Structural Mechanics}, Vol. 1. Dorset, England (Oct. 12-17, 1975):J.1--J.51,
  1975.

\bibitem{MR851383}
Vivette Girault and Pierre-Arnaud Raviart.
\newblock {\em Finite element methods for {N}avier-{S}tokes equations},
  volume~5 of {\em Springer Series in Computational Mathematics}.
\newblock Springer-Verlag, Berlin, 1986.
\newblock Theory and algorithms.

\bibitem{GG}
J.~Gopalakrishnan and J.~Guzm{\'a}n.
\newblock A second elasticity element using the matrix bubble.
\newblock {\em IMA J. Numer. Anal.}, 32(1):352--372, 2012.

\bibitem{MR3149075}
Johnny Guzm\'{a}n and Michael Neilan.
\newblock Symmetric and conforming mixed finite elements for plane elasticity
  using rational bubble functions.
\newblock {\em Numer. Math.}, 126(1):153--171, 2014.

\bibitem{HSV}
Antti Hannukainen, Rolf Stenberg, and Martin Vohral\'{\i}k.
\newblock A unified framework for a posteriori error estimation for the
  {S}tokes problem.
\newblock {\em Numer. Math.}, 122(4):725--769, 2012.

\bibitem{JM}
C.~Johnson and B.~Mercier.
\newblock Some equilibrium finite element methods for two-dimensional
  elasticity problems.
\newblock {\em Numer. Math.}, 30(1):103--116, 1978.

\bibitem{kim}
Kwang-Yeon Kim.
\newblock Guaranteed a posteriori error estimator for mixed finite element
  methods of linear elasticity with weak stress symmetry.
\newblock {\em SIAM J. Numer. Anal.}, 49(6):2364--2385, 2011.

\bibitem{LS1}
Philip~L. Lederer and Rolf Stenberg.
\newblock Energy norm analysis of exactly symmetric mixed finite elements for
  linear elasticity, 2021.

\bibitem{LV}
Marco Lonsing and R\"{u}diger Verf\"{u}rth.
\newblock A posteriori error estimators for mixed finite element methods in
  linear elasticity.
\newblock {\em Numer. Math.}, 97(4):757--778, 2004.

\bibitem{NecasH}
Jind\v{r}ich Ne\v{c}as and Ivan Hlav\'{a}\v{c}ek.
\newblock {\em Mathematical theory of elastic and elasto-plastic bodies: an
  introduction}, volume~3 of {\em Studies in Applied Mechanics}.
\newblock Elsevier Scientific Publishing Co., Amsterdam-New York, 1980.

\bibitem{MR25902}
W.~Prager and J.~L. Synge.
\newblock Approximations in elasticity based on the concept of function space.
\newblock {\em Quart. Appl. Math.}, 5:241--269, 1947.

\bibitem{netgen}
J.~Sch{\"o}berl.
\newblock {NETGEN An advancing front 2D/3D-mesh generator based on abstract
  rules}.
\newblock {\em Computing and Visualization in Science}, 1(1):41--52, 1997.

\bibitem{MR813691}
L.~R. Scott and M.~Vogelius.
\newblock Norm estimates for a maximal right inverse of the divergence operator
  in spaces of piecewise polynomials.
\newblock {\em RAIRO Mod\'{e}l. Math. Anal. Num\'{e}r.}, 19(1):111--143, 1985.

\bibitem{family}
Rolf Stenberg.
\newblock A family of mixed finite elements for the elasticity problem.
\newblock {\em Numer. Math.}, 53(5):513--538, 1988.

\bibitem{PPNM}
Rolf Stenberg.
\newblock Some new families of finite elements for the {S}tokes equations.
\newblock {\em Numer. Math.}, 56(8):827--838, 1990.

\bibitem{PPRAIRO}
Rolf Stenberg.
\newblock Postprocessing schemes for some mixed finite elements.
\newblock {\em RAIRO Mod\'el. Math. Anal. Num\'er.}, 25(1):151--167, 1991.

\bibitem{szabo1991finite}
B.~A. Szab{\'o} and I.~Babu{\v{s}}ka.
\newblock {\em Finite Element Analysis}.
\newblock A Wiley-Interscience publication. Wiley, 1991.

\bibitem{MR3059294}
R\"{u}diger Verf\"{u}rth.
\newblock {\em A posteriori error estimation techniques for finite element
  methods}.
\newblock Numerical Mathematics and Scientific Computation. Oxford University
  Press, Oxford, 2013.

\end{thebibliography}
\bibliographystyle{plain}

  \end{document}